\documentclass[12pt]{article}%
\usepackage{amsfonts}
\usepackage{mathrsfs}
\usepackage{fancyhdr}
\usepackage{comment}
\usepackage[a4paper, top=2.5cm, bottom=2.5cm, left=2.2cm, right=2.2cm]%
{geometry}
\usepackage{times}
\usepackage{amsmath}
\usepackage{changepage}
\usepackage{mathtools}
\usepackage{amssymb}
\usepackage{graphicx}%
\usepackage{amsthm}
\usepackage{tikz}
\usepackage{verbatim}
\usepackage[normalem]{ulem}
\usepackage[unicode=true]{hyperref}
\usepackage{cleveref}

\usetikzlibrary{arrows.meta}
\usetikzlibrary{positioning}
\newtheorem{theorem}{Theorem}[section]

\newtheorem{lemma}[theorem]{Lemma}

\newtheorem{definition}[theorem]{Definition}
\newtheorem{example}[theorem]{Example}

\newtheorem{remark}[theorem]{Remark}
\newtheorem{assumptionA}{A}
\newtheorem{assumptionB}{B}
\newtheorem{assumptionC}{C}
\newtheorem{assumptionD}{D}

\newcommand{\R}{\mathbb{R}}

\newcommand{\T}{\mathbb{T}}

\newcommand{\abs}[1]{\left \lvert #1 \right \rvert}

\newcommand{\bb}[1]{\mathbb{#1}}

\title{Some Remarks on the Asymptotic Behavior of Mean Field Games of Controls}
\author{Jameson Graber, Elizabeth Matter, Kyle Rosengartner}
\date{October 2026}

\begin{document}

\maketitle

\begin{abstract}
    We analyze the asymptotic behavior of several mean field games of controls. First, in the time-dependent case, we examine the behavior of solutions to discounted mean field games of controls as the discount factor $\lambda$ tends to infinity. Next, we study three stationary models of mean field games of controls, finding that solutions to the first-order discounted problem and the second-order ergodic problem can be obtained from solutions to the second-order discounted problem by taking $\nu \to 0$ and $\lambda \to 0$, respectively. Finally, we investigate the long-time average behavior of solutions to the time-dependent system with no discount, finding that solutions converge to the solution to the second-order ergodic problem as $T \to \infty$.
\end{abstract}

%-----%-----%-----%
%-----%-----%-----%
%-----%-----%-----%
%-----%-----%-----%
\section{Introduction}
%-----%-----%-----%
%-----%-----%-----%
%-----%-----%-----%
%-----%-----%-----%

The theory of mean field games of controls is closely related to the theory of mean field games, which was introduced independently by Lasry and Lions in \cite{lasry2007mean} and by Caines, Huang, and Malham\'e in \cite{huang2006large}. Mean field games (MFGs) and mean field games of controls (MFGCs) are both differential games in which players seek to minimize costs that depend on the collective behavior of all the players, and it is well-known that Nash equilibria can be characterized by solutions to coupled systems of PDE. The difference is that in a mean field game, this cost depends only on the distribution $m$ of players' states, whereas in a mean field game of controls, this cost depends on the joint distribution $\mu$ of players' states and controls. Thus, while the system of PDE characterizing a Nash equilibrium for a MFG consists only of a Hamilton-Jacobi equation satisfied by the optimized value function $u$ and a Fokker-Planck equation satisfied by the distribution of player states $m$, in the system corresponding to a MFGC, there is an additional fixed-point relation satisfied by the joint distribution $\mu$ of states and controls. This variant of a mean field game was initially called an \textit{extended mean field game} (see \cite{gomes2014existence,gomes2016extended}), but the terminology ``mean field game of controls" is now standard, cf.~\cite{cardaliaguet2018mean}.

Many papers have examined the asymptotic behavior of solutions to MFGs. For example, in \cite[Section 3]{bertucci2019some}, the authors give some conditions under which MFG models can approach agent based models (ABMs) as the discount factor $\lambda$ becomes very large. This inspired the authors of \cite{bardi2021convergence} to examine the asymptotic behavior of a rescaled MFG system as $\lambda \to \infty$ and $\nu \downarrow \nu_\infty \geq 0$, and a recent paper (see \cite{cirant2025uniqueness}) obtained quantitative convergence rates in the case that $\nu$ is fixed. See also \cite{degond2017mean} for related work. In \cite{ciampa2026rate,tang2025convergence}, the authors not only prove convergence as $\nu \to 0$ but also analyze the precise rates of convergence. The papers \cite{cardaliaguet2012long,cardaliaguet2013long,cirant2024long} studied the long-time average of solutions to second-order MFG systems as the time horizon $T \to \infty$ for both local and nonlocal couplings, the long-time average of solutions to first-order MFG systems as $T \to \infty$ was studied in \cite{bardi2024long,cannarsa2023weak,cardaliaguet2013KAM,cardaliaguet2015mean}, and the long-time average behavior of solutions to the master equation was studied in \cite{cardaliaguet2019long}. Finally, in the stationary setting, the papers \cite{gomes2019selection,iturriaga2025discretization,terai2026remarks} have investigated vanishing discount problems in first-order stationary MFG systems.

Compared to traditional MFGs, MFGCs have received far less attention in the literature, but they have gained attention in recent years, with most of the results focusing on well-posedness. In 2022, Kobeissi proved the well-posedness of second-order MFGCs in the periodic setting and the whole space under both monotone and non-monotone couplings (see \cite{kobeissi2022mean,kobeissi2022classical}). More recent papers have investigated the well-posedness of MFGCs under various boundary conditions, invariance and state constraints, and even with fractional diffusion replacing standard diffusion (see \cite{bo2025mean,bongini2024mean,graber2021note,graber2026mean,graber2025mean,graber2026ergodic}). The associated master equation was also recently studied in \cite{jackson2026unconditional,li2026master,mou2022mfgc}. Other noteworthy papers on special cases of MFGCs include \cite{camilli2025learning,chan2015bertrand,chan2017fracking,graber2024new,graber2018existence,graber2018variational,graber2023master,graber2020mean}. In this paper, we examine the asymptotic behavior of time-dependent MFGC systems as $\lambda \to \infty$, the convergence of second-order stationary MFGCs as the discount factor or viscosity vanishes, and the long-time average behavior of zero-discount MFGC systems as $T \to \infty$.

%-----%-----%-----%
%-----%-----%-----%
%-----%-----%-----%
%-----%-----%-----%
\subsection{Notation \& Preliminaries}
%-----%-----%-----%
%-----%-----%-----%
%-----%-----%-----%
%-----%-----%-----%

For $R > 0$, we will use $\mathcal{P}_{\infty,R}(\T^d \times \R^d)$ to denote the set of Borel probability measures $\mu \in \mathcal{P}(\T^d \times \R^d)$ with $\operatorname{supp} \mu \subseteq \{(x,y) \in \T^d \times \R^d : |y| \leq R\}$, and unless otherwise specified, each probability space will be endowed with the weak-* topology. For $\mathcal{P}_{\infty,R}(\T^d \times \R^d)$, this is induced by the Wasserstein metric
\begin{equation}
	W_1(\mu,\nu) = \inf\left\{\int \abs{x-y} d\pi(x,y) : \pi \in \Pi(\mu,\nu)\right\}
\end{equation}
where $\Pi(\mu,\nu)$ denotes the set of all couplings between $\mu$ and $\nu$.
By Kantorovitch duality, we also have the characterization
\begin{equation}
	W_1(\mu,\nu) = \sup\left\{\int \varphi d(\mu-\nu) : \|\varphi\|_{Lip}  \leq 1\right\}.
\end{equation}
Moreover, as in \cite{kobeissi2022mean,kobeissi2022classical}, we will define the quantity
\begin{equation} \label{Eq: Lambda-q}
    \Lambda_r(\mu) \coloneqq
    \begin{cases}
        \left(\int_{\Omega \times \R^d} |\alpha|^r d\mu(x,\alpha)\right)^\frac{1}{r}, &1 < r < \infty \\
        \sup\{|\alpha| : (x,\alpha) \in \text{supp} \mu\}, &r = \infty
    \end{cases}
\end{equation}
which will allow us to quantify the dependence of $H$ and its derivatives on $\mu$.

We will use the following function spaces: for non-negative integers $j,k$ we denote by $C^{j,k}$ the space of all functions $u(t,x)$ on $[0,T] \times \T^d$ that are $j$ times continuously differentiable with respect to $t$ and $k$ times continuously differentiable with respect to $x$.
For a fraction $\alpha \in (0,1)$, we denote by $C^{1+\alpha/2,2+\alpha}$ the usual parabolic H\"older space as introduced for instance in \cite{ladyzhenskaia1968linear}.
As for Sobolev spaces, we use similar notation as in \cite{ladyzhenskaia1968linear}, in particular $W^{1,2}_r$ denotes the space of all functions with weak derivatives up to order 1 in time and 2 in space, whose weak derivatives are all $L^r$ summable.
We similarly define the spaces $C^k$, $C^{k+\alpha}$, and $W^k_r$ for functions of $x$ alone.

If a function $\phi$ is semiconcave, i.e., $D_{xx}^2 \phi \leq C$, we will denote by $\operatorname{SC}(\phi)$ the smallest constant $C$ for which this inequality holds.

\begin{definition}
    For the MFGC system
    \begin{equation} \label{Eq: MFG 1}
        \begin{cases}
            -\partial_t u_\lambda -\nu\Delta u_\lambda + H(t,x,D_xu_\lambda,\mu_\lambda) + \lambda u_\lambda = f(t,x,m_\lambda(t)), &(t,x) \in (0,T) \times \T^d \\
            \partial_t m_\lambda - \nu\Delta m_\lambda - \operatorname{div} (m_\lambda D_pH(t,x,D_x u_\lambda,\mu_\lambda)) = 0, &(t,x) \in (0,T) \times \T^d \\
            \mu_\lambda = (I,-D_pH(t,\cdot,D_xu_\lambda(t,\cdot),\mu_\lambda)) \# m_\lambda, &t \in (0,T) \\
            m_\lambda(0,x) = m_0(x), \qquad u_\lambda(T,x) = g(x,m_\lambda(T)), &x \in \T^d
        \end{cases}
    \end{equation}
    with $\lambda \geq 0$, we distinguish between two different notions of solution.
    \begin{itemize}
        \item We will say $(u_\lambda,m_\lambda,\mu_\lambda) \in C^{0,1}([0,T] \times \T^d) \times C^0([0,T] \times \T^d) \times C^0(0,T;\mathcal{P}(\T^d \times \R^d))$ is a solution \textbf{in the sense of distributions} if $\mu$ satisfies the fixed-point relation for all $t \in [0,T]$ and the functions $u_\lambda$ and $m_\lambda$, respectively, solve the Hamilton-Jacobi and the Fokker-Planck equations in the sense of distributions.
        \item We will say that $(u_\lambda,m_\lambda,\mu_\lambda) \in C^{1,2}([0,T] \times \T^d) \times C^0([0,T] \times \T^d) \times C^0(0,T;\mathcal{P}(\T^d \times \R^d))$ is a \textbf{classical solution} if $(u_\lambda,m_\lambda,\mu_\lambda)$ is a solution in the sense of distributions and $u_\lambda$ satisfies the HJ equation in the classical sense.
    \end{itemize}
\end{definition}

Aside from these preliminaries, we specify that the constant $C$ appearing in many results denotes a generic constant that may change from line to line but depends only on the data, i.e., on the constants found in the assumptions (see Section \ref{Sec: Assumptions}).

%-----%-----%-----%
%-----%-----%-----%
%-----%-----%-----%
%-----%-----%-----%
\subsection{Summary of Our Results}
%-----%-----%-----%
%-----%-----%-----%
%-----%-----%-----%
%-----%-----%-----%

In this article, we will consider the asymptotic behavior of several models of MFGCs. We begin in Section \ref{Sec: Discounted MFGC} by considering the asymptotic behavior of two time-dependent MFGC systems as $\lambda \to \infty$ and $\nu_\lambda \to \nu_\infty \geq 0$, generalizing results from \cite{bardi2021convergence,bertucci2019some,cirant2025uniqueness} to the setting of mean field games of controls. In particular, we study the standard MFGC system \eqref{Eq: MFG 1} and the rescaled system
\begin{equation} \label{Eq: Rescaled MFG}
    \begin{cases}
        -\partial_t u_\lambda -\nu_\lambda\Delta u_\lambda + \lambda^{-1}H(t,x,\lambda D_xu_\lambda,\mu_\lambda) + \lambda u_\lambda = f(t,x,m_\lambda(t)), &(t,x) \in (0,\infty) \times \T^d \\
        \partial_t m_\lambda - \nu_\lambda\Delta m_\lambda - \operatorname{div} (m_\lambda D_pH(t,x,\lambda D_x u_\lambda,\mu_\lambda)) = 0, &(t,x) \in (0,\infty) \times \T^d \\
        \mu_\lambda = (I,-D_pH(t,\cdot,\lambda D_xu_\lambda(t,\cdot),\mu_\lambda)) \# m_\lambda, &t \in (0,\infty) \\
        m_\lambda(0,x) = m_0(x), &x \in \T^d
    \end{cases}
\end{equation}
See Sections \ref{Sec: ABM} and \ref{Sec: Aggregation}, respectively. In the former case, we find that $\|u_\lambda(t,\cdot)\|_{C^1(\T^d)}$ vanishes locally uniformly in $[0,T)$ as $\lambda \to \infty$, and in the latter case, we find that $(\lambda u_\lambda)_\lambda$ is uniformly bounded in $C^1(\T^d)$ and uniformly semiconcave. In both cases, we find that the sequence $(m_\lambda)_\lambda$ is relatively compact, and more importantly, we show that any limit point of \eqref{Eq: MFG 1} solves
\begin{equation} \label{Eq: ABM}
    \begin{cases}
        \partial_t m - \nu_\infty\Delta m - \operatorname{div} (mD_pH(t,x,0,\mu)) = 0, &(t,x) \in (0,T) \times \T^d \\
        \mu = (I,-D_pH(t,\cdot,0,\mu)) \# m, &t \in (0,T) \\
        m(0,x) = m_0(x), &x \in \T^d
    \end{cases}
\end{equation}
and any limit point of \eqref{Eq: Rescaled MFG} solves
\begin{equation} \label{Eq: Aggregation}
    \begin{cases}
        \partial_t m - \nu_\infty\Delta m - \operatorname{div} (mD_pH(t,x,D_xf(t,x,m),\mu)) = 0, &(t,x) \in (0,\infty) \times \T^d \\
        \mu = (I,-D_pH(t,\cdot,D_xf(t,\cdot,m),\mu)) \# m, &t \in (0,\infty) \\
        m(0,x) = m_0(x), &x \in \T^d
    \end{cases}
\end{equation}
Additionally, we observe that under additional assumptions, we can obtain uniqueness of the limits (and hence full convergence of each sequence) as well as quantitative convergence rates in the case where $\nu_\lambda \equiv \nu_0$ for some $\nu_0 > 0$, using methods found in \cite{cirant2025uniqueness}.

Crucial to our arguments in Sections \ref{Sec: ABM} and \ref{Sec: Aggregation} are the estimates obtained in \ref{Sec: Estimates}. In addition to our convergence results, which are the purpose of this article, these estimates allow us to prove the existence of solutions to MFGCs on infinite time horizons in Section \ref{Sec: Existence}.

In Section \ref{Sec: Stationary}, we show the well-posedness of the stationary systems
\begin{equation} \label{Eq: Stationary Discounted MFGC}
    \begin{cases}
        \lambda u - \nu\Delta u + H(x,D_xu,\mu) = f(x,m), &x \in \T^d \\
        \lambda m - \nu\Delta m - \operatorname{div} (mD_pH(x,D_xu,\mu)) = \lambda, &x \in \T^d \\
        \mu = (I,-D_pH(\cdot,D_xu,\mu)) \# m
    \end{cases}
\end{equation}
\begin{equation} \label{Eq: Ergodic MFGC}
    \begin{cases}
        -\nu\Delta u + H(x,D_xu,\mu) + \rho = f(x,m), &x \in \T^d \\
        \nu\Delta m + \operatorname{div} (mD_pH(x,D_xu,\mu)) = 0, &x \in \T^d \\
        \mu = (I,-D_pH(\cdot,D_xu,\mu)) \# m \\
        m \geq 0, \qquad \int_{\T^d} m \, dx = 1
    \end{cases}
\end{equation}
\begin{equation} \label{Eq: Stationary First-Order MFGC}
    \begin{cases}
        \lambda u +H(x,D_xu,\mu) = f(x,m), &x \in \T^d \\
        \lambda m - \operatorname{div} (mD_pH(x,D_xu,\mu)) = \lambda, &x \in \T^d \\
        \mu = (I,-D_pH(\cdot,D_xu,\mu)) \# m
    \end{cases}
\end{equation}
and show that the solutions to \eqref{Eq: Ergodic MFGC}, \eqref{Eq: Stationary First-Order MFGC} are obtained from those of \eqref{Eq: Stationary Discounted MFGC} by taking the limit as $\lambda \to 0$ or $\nu \to 0$, respectively. This is summarized by Figure \ref{Fig: Stationary Problems}. In addition, we investigate the rate of convergence for the vanishing viscosity problem $\eqref{Eq: Stationary Discounted MFGC} \to \eqref{Eq: Stationary First-Order MFGC}$. In particular, we find that
$$
\begin{aligned}
    \int_{\T^d} (m_{\lambda,\nu} + m_\lambda)|D_pH(x,D_xu_{\lambda,\nu},\mu_{\lambda,\nu}) - D_pH(x,D_xu_\lambda,\mu_\lambda)|^2 dx \\
    + \int_{\T^d} (f(x,m_{\lambda,\nu}) - f(x,m_\lambda))d(m_{\lambda,\nu} - m_\lambda) \leq C\sqrt{\nu}.
\end{aligned}
$$

\begin{figure}[htbp]
    \centering
    \begin{tikzpicture}[
    box/.style={draw, rectangle, minimum width=2.4cm, minimum height=1.2cm, align=center},
    arrow/.style={->, thick}
]

% Top row
\node[box] (A) {$
    \begin{cases}
        \lambda u_{\lambda,\nu} - \nu\Delta u_{\lambda,\nu} + H(x,D_xu_{\lambda,\nu},\mu_{\lambda,\nu}) = f(x,m_{\lambda,\nu}) \\
        \lambda m_{\lambda,\nu} - \nu\Delta m_{\lambda,\nu} - \operatorname{div} (m_{\lambda,\nu}D_pH(x,D_xu_{\lambda,\nu},\mu_{\lambda,\nu})) = \lambda \\
        \mu_{\lambda,\nu} = (I,-D_pH(\cdot,D_xu_{\lambda,\nu},\mu_{\lambda,\nu})) \# m_{\lambda,\nu}
    \end{cases}
    $
};

% Bottom row
\node[box, below left=1cm and -5.4cm of A] (B) {$
     \begin{cases}
        - \nu \Delta u_\nu + H(x,Du_\nu,\mu_\nu) + \rho_\nu = f(x, m_\nu), \\
        \nu\Delta m_\nu + \operatorname{div} (m_\nu D_pH(x,Du_\nu,\mu_\nu)) = 0 \\
        \mu_\nu = (I,-D_pH(\cdot,D_xu_\nu,\mu_\nu)) \# m_\nu
    \end{cases}
    $
    };
    
\node[box, below right=1cm and -5.4cm of A] (C) {
    $
     \begin{cases}
        \lambda u_\lambda + H(x,Du_\lambda,\mu_\lambda) = f(x, m_\lambda), \\
        \lambda m_\lambda - \operatorname{div} (m_\lambda D_pH(x,Du_\lambda,\mu_\lambda)) = \lambda \\
        \mu_\lambda = (I,-D_pH(\cdot,D_xu_\lambda,\mu_\lambda)) \# m_\lambda
    \end{cases}
    $
};

\draw[arrow] (A) -- node[left=0.5] {$\lambda \to 0$} (B);

% Vertical arrows
\draw[arrow] (A) -- node[right=0.5cm] {$\nu \to 0$} (C);

\end{tikzpicture}
\caption{Summary of Convergence for Stationary Systems}
\label{Fig: Stationary Problems}
\end{figure}

Finally, in the spirit of \cite{cardaliaguet2012long,cardaliaguet2013long}, Section \ref{Sec: Long-Time Average} is devoted to investigating the long-time average behavior of solutions to the zero-discount MFGC system
\begin{equation} \label{Eq: MFGC with No Discount Factor}
    \begin{cases}
        -\partial_t u -\nu\Delta u + H(t,x,D_xu,\mu) = f(t,x,m(t)), &(t,x) \in (0,T) \times \T^d \\
        \partial_t m - \nu\Delta m - \operatorname{div} (mD_pH(t,x,D_x u,\mu)) = 0, &(t,x) \in (0,T) \times \T^d \\
        \mu = (I,-D_pH(t,\cdot,D_xu(t,\cdot),\mu)) \# m, &t \in (0,T) \\
        m(0,x) = m_0(x), \qquad u(T,x) = g(x,m(T)), &x \in \T^d
    \end{cases}
\end{equation}
as $T \to \infty$. In particular, we find that for a particular class of Hamiltonians, the long-time average of solutions to \eqref{Eq: MFGC with No Discount Factor} converges to \eqref{Eq: Ergodic MFGC}.

%-----%-----%-----%
%-----%-----%-----%
%-----%-----%-----%
%-----%-----%-----%
\subsection{Assumptions}
\label{Sec: Assumptions}
%-----%-----%-----%
%-----%-----%-----%
%-----%-----%-----%
%-----%-----%-----%

In this section, we establish the main assumptions that we will make in the later sections. We remark that the constants $C_0,q,\beta_0,L_1$ are fixed .

First, to ensure the existence of solutions to \eqref{Eq: MFG 1} and \eqref{Eq: Rescaled MFG} in the sense of distributions, we make the following assumptions (see \cite{kobeissi2022mean}).

\begin{assumptionA} \label{A: A1}
    $L: [0,T] \times \T^d \times \R^d \times \mathcal{P}(\T^d \times \R^d) \to \R$  is differentiable with respect to $(x,\alpha)$; $L$ and its derivatives are continuous on $[0,T] \times \T^d \times \R^d \times \mathcal{P}_{\infty,R}(\T^d \times \R^d)$ for any $R > 0$; we recall that $\mathcal{P}_{\infty,R}(\T^d \times \R^d)$ is endowed with the weak* topology on measures; we use the notation $D_xL$, $D_\alpha L$ and $D_{(x,\alpha)}L$ for respectively the first-order derivatives of L with respect to $x$, $\alpha$ and $(x,\alpha)$.
\end{assumptionA}
\begin{assumptionA} \label{A: A2}
    The maximum in 
    \begin{equation}
        \sup_\alpha\{-\alpha \cdot p - L(t,x,\alpha,\mu)\}
    \end{equation}
    is achieved at a unique $\alpha \in \R^d$.
\end{assumptionA}
By \ref{A: A2} we may define the \emph{Hamiltonian}
$$H(t,x,p,\mu)=\sup_{\alpha\in\bb{R}^d}\{-\alpha \cdot p - L(t,x,\alpha,\mu)\},$$ and the maximizer is $\alpha^*=-D_pH(t,x,p,\mu)$.
\begin{assumptionA} \label{A: A3}
    For any $\mu_1,\mu_2 \in \mathcal{P}(\T^d \times \R^d)$, $L$ satisfies
    $$\int_{\T^d \times \R^d} (L(t,x,\alpha,\mu_1) - L(t,x,\alpha,\mu_2))d(\mu_1 - \mu_2)(x,\alpha) \geq 0.$$
\end{assumptionA}
\begin{assumptionA} \label{A: A4}
    $L(t,x,\alpha,\mu) \geq C_0^{-1}|\alpha|^{q'} - C_0(1 + \Lambda_{q'}(\mu)^{q'})$ where $q' = \frac{q}{q-1}$ for some $q \in (1,\infty)$.
\end{assumptionA}
\begin{assumptionA} \label{A: A5}
    $|L(t,x,\alpha,\mu)| + |D_xL(t,x,\alpha,\mu)| \leq C_0(1 + |\alpha|^{q'} + \Lambda_{q'}(\mu)^{q'})$.
\end{assumptionA}
\begin{assumptionA} \label{A: A6}
    The maps $g: \mathcal{P}(\T^d) \to C^{2+\beta_0}(\T^d)$ is continuous, $f: \mathcal{P}(\T^d) \to C^{1/2,1}([0,T] \times \T^d)$ is H\"{o}lder continuous (uniformly in $t$), and
    $$\|m_0\|_{C^{\beta_0}} + \underset{m}{\sup} \|f(\cdot,\cdot,m)\|_{C^{1/2,1}} + \underset{m}{\sup} \|g(\cdot,m)\|_{C^{2+\beta_0}} \leq C_0,$$
    where the supremums are taken over $t \in [0,T]$ and $m \in \mathcal{P}(\T^d)$.
\end{assumptionA}

To obtain classical solutions and prove the existence of solutions on infinite time horizons, we will use the following additional assumptions.

\begin{assumptionB} \label{A: DpH is Lipschitz in alpha}
    For every $R > 0$, $p \in \R^d$, and $\mu_1,\mu_2 \in \mathcal{P}_{\infty,R}(\T^d \times \R^d)$ with $\mu_i = (I,\alpha^{\mu_i}) \# \rho$, we have
    $$|D_pH(t,x,p,\mu_1) - D_pH(t,x,p,\mu_2)| \leq L_1\|\alpha^{\mu_1} - \alpha^{\mu_2}\|_{L^{q}(\rho)}$$
    for some $L_1 \in (0,1)$.
\end{assumptionB}
\begin{assumptionB} \label{A: DpH is Lipschitz in m}
    Fix $R > 0$, $p \in \R^d$ with $|p| \leq R$, and $\mu_1,\mu_2 \in \mathcal{P}_{\infty,R}(\T^d \times \R^d)$ with $\mu_i = (I,\alpha) \# \rho_i$ where $\|\alpha\|_{L^\infty(\T^d)} \leq R$. There exists $C_R > 0$ depending only on $R$ so that if $\rho_1,\rho_2 \in L^\infty(\T^d)$, then
    $$|D_pH(t,x,p,\mu_1) - D_pH(t,x,p,\mu_2)| \leq C_R\|\rho_1 - \rho_2\|_{L^\infty(\T^d)}$$
    and if $\alpha \in W^{1,\infty}(\T^d)$, then
    $$|D_pH(t,x,p,\mu_1) - D_pH(t,x,p,\mu_2)| \leq C_R(\|\alpha\|_{W^{1,\infty}(\T^d)} + 1)W_1(\rho_1,\rho_2).
    $$
\end{assumptionB}
\begin{assumptionB} \label{A: Regularity of DpH}
    Fix $R > 0$, $p_1,p_2 \in \R^d$ with $|p_1| + |p_2| \leq R$, and $\mu \in \mathcal{P}_{\infty,R}(\T^d \times \R^d)$. There exists $C_R > 0$ depending only on $R$ such that
    $$
    \begin{aligned}
        |D_pH(t,x,p_1,\mu) - D_pH(s,y,p_2,\mu)| \leq C_R(|p_1 - p_2| + |x-y| + |t-s|)
    \end{aligned}
    $$
\end{assumptionB}
\begin{assumptionB} \label{A: Regularity of H}
    Fix $R > 0$, $p \in \R^d$ with $|p| \leq R$, and $\mu_1,\mu_2 \in \mathcal{P}_{\infty,R}(\T^d \times \R^d)$. There exists $C_R > 0$ depending only on $R$ such that
    $$|H(t,x,p,\mu_1) - H(s,x,p,\mu_2)| \leq C_R(|t-s| + W_1(\mu_1,\mu_2)).$$
\end{assumptionB}

For the rescaled system, we will also need to make the following assumptions to obtain the desired convergence results.

\begin{assumptionC} \label{A: f is semiconcave}
    $f(t,x+h,m) + f(t,x-h,m) - 2f(t,x,m) \leq C_0|h|^2$ for all $t \in (0,T)$, $m \in \mathcal{P}(\T^d)$.
\end{assumptionC}
\begin{assumptionC} \label{A: Second Derivatives of H}
    $H$ is twice differentiable with respect to $(x,p)$, and the functions $D_{xx}^2H$, $D_{px}^2H$, and $D_{pp}^2H$ are continuous on $[0,T] \times \T^d \times \R^d \times \mathcal{P}(\T^d \times \R^d)$. Moreover, for all $x,h \in \T^d$, $t \in (0,T)$, $p \in \R^d$ with $|p| \leq R$, and $\mu \in \mathcal{P}_{\infty,R}(\T^d \times \R^d)$, we have
    \begin{enumerate}
        \item $H(t,x+h,p,\mu) + H(t,x-h,p,\mu) - 2H(t,x,p,\mu) \geq -C_R|h|^2
        $
        \item $D_{pp}^2H(t,x,p,\mu) \geq \frac{1}{C_R}$
        \item $|D_{pp}^2H(t,x,p,\mu)| \leq C_R$.
    \end{enumerate}
\end{assumptionC}

Finally, in our analysis of stationary systems and the long-time average behavior of solutions to \eqref{Eq: MFGC with No Discount Factor}, we will assume the following:

\begin{assumptionD} \label{A: Bounds for H}
    The functions $H$, $L$, and $f$ do not depend on $t$. Furthermore, we have
    \begin{equation*}\label{Eq: Growth of H}
        \lim_{|p|\to \infty} \inf_{x\in \T^d} \left(\frac{1}{2}H(x,p,\mu)^2 + D_xH(x,p,\mu) \cdot p\right) = \infty,  \qquad \lim_{|p|\to \infty} \inf_{x \in \T^d} \frac{H(x,p,\mu)^2}{|p|} = \infty
    \end{equation*}
    for each $\mu$, and
    $$H(x,p,\mu) \geq \frac{1}{C_0}|p|^q - C_0(1 + \Lambda_{q'}(\mu)^{q'}).$$
\end{assumptionD}

\begin{assumptionD} \label{A: DpxH = 0 or H is uniformly convex}
    $H$ is twice differentiable with respect to $(x,p)$, and its derivatives $D_{xx}^2H$, $D_{pp}^2H$, and $D_{px}^2H$ are continuous on $\T^d \times \R^d \times \mathcal{P}(\T^d \times \R^d)$. Furthermore, $H$ satisfies
    $$|D_{pp}^2H(x,p,\mu)| + |D_{px}^2H(x,p,\mu)| + |D_{xx}^2H(x,p,\mu)| \leq C_R$$
    and
    $$H(x,p_1,\mu) - H(x,p_2,\mu) - D_pH(x,p_2,\mu) \cdot (p_1 - p_2) \geq \frac{1}{C_R}|p_1 - p_2|^2$$
    for all $x \in \T^d$, $p,p_1,p_2 \in \R^d$, and $\mu \in \mathcal{P}_{\infty,R}(\T^d \times \R^d)$.
\end{assumptionD}

\begin{assumptionD} \label{A: U}
    For all $m_1,m_2 \in \mathcal{P}(\T^d)$, we have
    $$\int_{\T^d} (f(x,m_1) - f(x,m_2))d(m_1 - m_2) \geq 0.$$
\end{assumptionD}

\begin{assumptionD} \label{A: L is Uniformly Convex}
    For all $x \in \T^d$, $\alpha_1,\alpha_2 \in \R^d$, and $\mu \in \mathcal{P}_{\infty,R}(\T^d \times \R^d)$
    $$L(x,\alpha_1,\mu) - L(x,\alpha_2,\mu) - D_\alpha L(x,\alpha_2,\mu) \cdot (\alpha_1 - \alpha_2) \geq \frac{1}{C_R}|\alpha_1 - \alpha_2|^2.$$
\end{assumptionD}

\begin{assumptionD} \label{A: Lipschitz Assumption on DpH Inverse}
    For all $p_1,p_2 \in \R^d$ with $|p_1| + |p_2| \leq R$ and $\mu \in \mathcal{P}_{\infty,R}(\T^d \times \R^d)$, we have
    $$|D_pH(p_1,\mu) - D_pH(p_2,\mu)| \geq \frac{1}{C_R}|p_1 - p_2|.$$
\end{assumptionD}

\begin{assumptionD} \label{A: No x-dependence}
    The functions $H$, $L$, and $f$ do not depend on $x$.
\end{assumptionD}

\begin{assumptionD} \label{A: Example Hamiltonian}
    The coupling terms $f,g$ satisfy $\underset{m}{\sup}(\|f(\cdot,m)\|_{W^2_\infty(\T^d)} + \|g(\cdot,m)\|_{C^{3+\beta}(\T^d)}) \leq C_0$, and our Hamiltonian takes the form $H(p,\mu) = \psi(\mu)|p + \varphi(\mu)|^2 + V(\mu)$, where
    \begin{enumerate}
        \item $\psi,\varphi,V$ are continuous
        \item $\frac{1}{C_0} \leq \psi \leq C_0$
        \item $|\varphi(\mu)|^2 + |V(\mu)| \leq C_0(1 + \Lambda_2(\mu)^2)$
    \end{enumerate}
\end{assumptionD}

%-----%-----%-----%
%-----%-----%-----%
%-----%-----%-----%
%-----%-----%-----%
\subsection{Motivating Examples}
%-----%-----%-----%
%-----%-----%-----%
%-----%-----%-----%
%-----%-----%-----%

To motivate our analysis, in this section, we discuss some examples of Hamiltonians that satisfy the assumptions made in the previous section. In particular, we observe that Hamiltonian-Lagrangian pairs of the form
\begin{equation} \label{Eq: Abstract Example}
    \begin{cases}
        H(t,x,p,\mu) = \psi(t,x,\mu)|p + \varphi(t,x,\mu)|^q + V(t,x,\mu) \\
        L(t,x,\alpha,\mu) = \left(q^{1-q'} - q^{-q'}\right)\psi(t,x,\mu)^{1-q'}|\alpha|^{q'} + \alpha \cdot \varphi(t,x,\mu) - V(t,x,\mu)
    \end{cases}
\end{equation}
satisfy these requirements under appropriate assumptions on $q$ and the functions $\psi$, $\varphi$, and $V$.

\begin{example} \label{Ex: Linear Demand Model}
    The first and most basic example we will consider is that of the linear demand system introduced in \cite{chan2015bertrand}, for which the associated Hamiltonian-Lagrangian pair is given by
    \begin{equation}
        \begin{cases}
            H^{lin}(p,\mu) = \frac{1}{4}\left(p + \frac{\varepsilon}{1 + \varepsilon}\overline{\alpha} - \frac{1}{1 + \varepsilon}\right)^2 \\
            L^{lin}(\alpha,\mu) = \alpha^2 + \frac{\varepsilon}{1 + \varepsilon}\alpha\overline{\alpha} - \frac{1}{1 + \varepsilon}\alpha
        \end{cases}
    \end{equation}
    where $\overline{\alpha} = \int_{\T \times \R} \widetilde{\alpha} d\mu(y,\widetilde{\alpha})$ (see also \cite{graber2018existence,graber2018variational}) and $\varepsilon>0$. It is straightforward to check that all of the assumptions on $H$ and $L$ are satisfied by this pair. 
\end{example}
For simplicity, this example is only in $d=1$. The results in Sections \ref{Sec: Stationary}-\ref{Sec: Long-Time Average}, which require $d\geq 2$ will not hold, however the following examples are for a general dimension $d.$
\begin{example} \label{Ex: Generalized Linear Demand Model}
    Generalizing the first example, we can consider
    \begin{equation}
        \begin{cases}
            H(t,x,p,\mu) = \phi_1(t,x)^{-1}\left|p + \frac{\varepsilon}{1 + \varepsilon}\overline{\alpha} + \phi_2(t,x)\right|^2 + V(t,x,\mu) \\
            L(t,x,\alpha,\mu) = \frac{1}{4}\phi_1(t,x)|\alpha|^2 + \alpha \cdot \left(\frac{\varepsilon}{1 + \varepsilon}\overline{\alpha} + \phi_2(t,x)\right) - V(t,x,\mu)
        \end{cases}
    \end{equation}
    where $\overline{\alpha} = \int_{\T^d \times \R^d} \widetilde{\alpha} d\mu(y,\widetilde{\alpha})$. If we have that $\phi_1,\phi_2,V$ are sufficiently smooth, $\phi_1 > \frac{2\varepsilon}{1+\varepsilon}$,
    $$\underset{t}{\sup} \|V(t,\cdot,\mu)\|_{C^1(\T^d)} \leq C(1 + \Lambda_{q'}(\mu)^{q'}), \qquad  \underset{t}{\sup}\left(\|\phi_1(t,\cdot)\|_{C^2(\T^d)} + \|\phi_2(t,\cdot)\|_{C^2(\T^d)}\right) \leq C,$$
    $V$ is Lipschitz in $W^1$, and 
   $$\int_{\T^d \times \R^d} \left(V(t,x,\mu_2) - V(t,x,\mu_1)\right) d(\mu_1 - \mu_2) \geq 0,$$
    then it is easy to confirm that such Hamiltonian-Lagrangian pairs satisfy A\ref{A: A1}-\ref{A: A5}, B\ref{A: DpH is Lipschitz in alpha}-\ref{A: Regularity of H}, and C\ref{A: Second Derivatives of H}.
\end{example}

\begin{example}
    We observe that for the example in \eqref{Eq: Abstract Example} to satisfy conditions A\ref{A: A1}-\ref{A: A2}, A\ref{A: A4}-\ref{A: A5}, and D\ref{A: Bounds for H}, we only need to assume that $\psi,\varphi,V$ are constant in $t$ and sufficiently smooth (e.g., differentiable in $x$ with $\psi,\varphi,V$ and their derivatives continuous on $\T^d \times \R^d \times \mathcal{P}(\T^d \times \R^d)$), that $\frac{1}{C} \leq \psi \leq C$ for some constant $C \geq 1$, and that
    $$\|V(\cdot,\mu)\|_{C^1(\T^d)} + \|\psi(\cdot,\mu)\|_{C^1(\T^d)} + \|\varphi(\cdot,\mu)\|_{C^1(\T^d)}^q \leq C(1 + \Lambda_{q'}(\mu)^{q'}).$$
    If $q = 2$ and we assume more regularity on $\psi,\varphi,V$, then this example also satisfies C\ref{A: Second Derivatives of H}, D\ref{A: DpxH = 0 or H is uniformly convex}, and D\ref{A: L is Uniformly Convex}-\ref{A: Lipschitz Assumption on DpH Inverse}. To satisfy condition A\ref{A: A3} is more complicated, but we point out a few sufficient conditions:
    \begin{enumerate}
        \item The functions $\varphi,V$ do not depend on $\mu$ and
        $$\int_{\T^d \times \R^d} \left(\psi(x,\mu_1)^{1-q'} - \psi(x,\mu_2)^{1-q'}\right)|\alpha|^{q'} d(\mu_1 - \mu_2) \geq 0.$$
        If $\psi$ does not depend on $x$, this reduces to $(\psi(\mu_1) - \psi(\mu_2)) \cdot (\Lambda_{q'}(\mu_1) - \Lambda_{q'}(\mu_2)) \leq 0$.
        \item The functions $\psi,V$ do not depend on $\mu$ and
        $$\int_{\T^d \times \R^d} \alpha \cdot \left(\varphi(x,\mu_1) - \varphi(x,\mu_2)\right) d(\mu_1 - \mu_2) \geq 0.$$
        If $\varphi$ does not depend on $x$, this reduces to $(\varphi(\mu_1) - \varphi(\mu_2)) \cdot \int_{\Omega \times \R^d} \widetilde{\alpha} d(\mu_1 - \mu_2)(x,\widetilde{\alpha}) \geq 0$.
        \item The functions $\psi,\varphi$ do not depend on $\mu$ and
        $$\int_{\T^d \times \R^d} \left(V(x,\mu_2) - V(x,\mu_1)\right) d(\mu_1 - \mu_2) \geq 0,$$
        which is trivial if $V$ does not depend on $x$.
    \end{enumerate}
\end{example}

%-----%-----%-----%
%-----%-----%-----%
%-----%-----%-----%
%-----%-----%-----%
\section{Convergence of Time Dependent MFGC with Large Discount}
\label{Sec: Discounted MFGC}
%-----%-----%-----%
%-----%-----%-----%
%-----%-----%-----%
%-----%-----%-----%

%-----%-----%-----%
%-----%-----%-----%
%-----%-----%-----%
%-----%-----%-----%
\subsection{A Priori Estimates for Discounted Systems}
\label{Sec: Estimates}
%-----%-----%-----%
%-----%-----%-----%
%-----%-----%-----%
%-----%-----%-----%

In this section, we adapt the arguments from \cite{kobeissi2022mean} to obtain a priori estimates for solutions to discounted mean field games of controls. We also introduce the parameter $\gamma$ to our calculations to allow us to obtain estimates for both \eqref{Eq: MFG 1} and the rescaled system \eqref{Eq: Rescaled MFG} simultaneously.

\begin{lemma} \label{Lem: General C1 Estimate}
    Assume A\ref{A: A1}-\ref{A: A6} hold. Fix $\nu > 0$, $\lambda > 0$, $\Omega \in \{\T^d,\R^d\}$, and suppose $(u_\lambda,m_\lambda,\mu_\lambda)$ solves
    \begin{equation} \label{Eq: General Discounted MFG}
        \begin{cases}
            -\partial_t u_\lambda -\nu\Delta u_\lambda + H(t,x,D_xu_\lambda,\mu_\lambda) + \lambda u_\lambda = \gamma f(t,x,m_\lambda(t)), &(t,x) \in (0,T) \times \Omega \\
            \partial_t m_\lambda - \nu\Delta m_\lambda - \operatorname{div} (m_\lambda D_pH(t,x,D_x u_\lambda,\mu_\lambda)) = 0, &(t,x) \in (0,T) \times \Omega \\
            \mu_\lambda = (I,-D_pH(t,\cdot,D_xu_\lambda(t,\cdot),\mu_\lambda)) \# m_\lambda, &t \in (0,T) \\
            m_\lambda(0,x) = m_0(x), \qquad u_\lambda(T,x) = \gamma g(x,m_\lambda(T)), &x \in \Omega
        \end{cases}
    \end{equation}
    where either $\gamma = 1$, $T \in (0,\infty]$ or $\gamma \in \{1,\lambda\}$, $T = \infty$. If A\ref{A: A1}-\ref{A: A6} hold (where $\T^d$ is replaced with $\Omega$), then there exists some $C > 0$ independent of $\nu,\lambda,T$ such that
    $$\|u_\lambda(t,\cdot)\|_\infty \leq
    \begin{cases}
        C(e^{\lambda(t-T)} + \lambda^{-1}), &\text{if } \gamma = 1, T < \infty \\
        C\lambda^{-1}, &\text{if } \gamma = 1, T = \infty \\
        C\max\{1,\lambda^{-1}\}, &\text{if } \gamma = \lambda, T = \infty
    \end{cases}$$
    and
    $$\|D_xu_\lambda(t,\cdot)\|_\infty \leq
    \begin{cases}
        C(e^{\frac{\lambda}{2}(t-T)} + \lambda^{-\frac{1}{2}}), &\text{if } \gamma = 1, T < \infty \\
        C\lambda^{-\frac{1}{2}}, &\text{if } \gamma = 1, T = \infty \\
        C\max\{1,\lambda^{-\frac{1}{2}}\}, &\text{if } \gamma = \lambda, T = \infty
    \end{cases}$$
\end{lemma}

\begin{proof}
    By the optimality of the control $\alpha_\lambda(t,\cdot) = -D_pH(t,\cdot,D_xu_\lambda(t,\cdot),\mu_\lambda(t))$, we get
    $$
    \begin{aligned}
        &\int_t^T e^{-\lambda s} \int_\Omega L(s,x,\alpha,\mu_\lambda) d(\mu_\lambda - \widetilde{\mu}) ds \\
        &\leq \gamma \int_t^T e^{-\lambda s} \int_\Omega f(s,x,m_\lambda(s)) d(\widetilde{m} - m_\lambda) ds + \gamma e^{-\lambda T}\int_\Omega g(x,m_\lambda(T)) d(\widetilde{m}(T) - m_\lambda(T)) \\
        &\leq 2\gamma C_0\lambda^{-1}(e^{-\lambda t} - e^{-\lambda T}) + 2C_0\gamma e^{-\lambda T}
    \end{aligned}
    $$
    where $\widetilde{m}(s) = \mathcal{L}(\widetilde{X}_t)$, $\widetilde{\mu}(s) = \widetilde{m}(s) \otimes \delta_0$, and
    $$d\widetilde{X}_s = \sqrt{2\nu}dB_s, \qquad \widetilde{X}_t \sim m(t).$$
    Next, by A\ref{A: A3}-\ref{A: A5}, we have
    $$
    \begin{aligned}
        \int_t^T e^{-\lambda s}\Lambda_{q'}(\mu_\lambda(s))^{q'} ds &= \int_t^T e^{-\lambda s}\int_\Omega |\alpha|^{q'} d\mu_\lambda(s) ds \\
        &\leq C_0\int_t^T e^{-\lambda s} \int_\Omega (L(s,x,\alpha,\widetilde{\mu}) + C_0) d\mu_\lambda(s) ds \\
        &\leq C_0\int_t^T e^{-\lambda s} \left(\int_\Omega L(s,x,\alpha,\mu_\lambda) d(\mu_\lambda - \widetilde{\mu}) + \int_\Omega L(s,x,\alpha,\widetilde{\mu}) d\widetilde{\mu}\right) ds \\
        &\qquad + C_0^2\lambda^{-1}(e^{-\lambda t} - e^{-\lambda T}) \\
        &\leq 2(1 + \gamma) C_0^2\lambda^{-1}(e^{-\lambda t} - e^{-\lambda T}) + 2C_0^2\gamma e^{-\lambda T}.
    \end{aligned}
    $$
    By the maximum principle, we get
    $$
    \begin{aligned}
        \|u_\lambda(t)\|_\infty &\leq \gamma e^{\lambda(t-T)}\|g(x,m_\lambda(T))\|_\infty \\
        &\qquad + e^{\lambda t}\int_t^T e^{-\lambda s}(\|H(s,x,0,\mu_\lambda)\|_\infty + \gamma\|f(s,x,m_\lambda(s))\|_\infty)ds \\
        &\leq (C_0 + 2C_0^3)\gamma e^{\lambda(t-T)} + (1+\gamma)(C_0 + 2C_0^3)\lambda^{-1}(1 - e^{\lambda(t-T)}).
    \end{aligned}
    $$
    If $T = \infty$, the first term is dropped and we use $\int_t^\infty e^{\lambda(t-s)} ds = \lambda^{-1}$, thus yielding our estimates on $u_\lambda$.
    
    For simplicity of presentation, for the gradient estimate, we will assume $u_\lambda$ is smooth (otherwise, we would take a convolution as in \cite{kobeissi2022mean}). Considering first the case where $\lambda \geq 1$, we define
    $$\varphi(v) = \exp(\exp(-1-R^{-1}v)) \qquad \text{for } |v| \leq R$$
    where $R \geq \|u_\lambda\|_\infty + (1+C_0)\underset{t}{\sup}\int_t^T e^{\lambda(t-s)}\Lambda_{q'}(\mu_\lambda(s))^{q'} ds$ to be chosen later, and
    $$w(t,x) = \varphi\left(u_\lambda(t,x) + (1+C_0)\int_t^T e^{\lambda(t-s)}\Lambda_{q'}(\mu_\lambda(s))^{q'} ds\right)|D_xu_\lambda(t,x)|^2.$$
    Through tedious computations, we get
    $$
    \begin{aligned}
        &-\partial_tw - \nu\Delta w + D_xw \cdot D_pH - 2\nu\frac{\varphi'}{\varphi}D_xw \cdot D_xu_\lambda + 2\nu\varphi|D_{xx}^2u_\lambda|^2 + 2\lambda w \\
        &= \frac{\varphi'}{\varphi}w\bigg[\gamma f - \lambda u_\lambda + D_xu_\lambda \cdot D_pH - H - \lambda(1+C_0)\int_t^T e^{\lambda(t-s)}\Lambda_{q'}(\mu_\lambda(s))^{q'} ds \\
        &\qquad\qquad + (1 + C_0)\Lambda_{q'}(\mu_\lambda)^{q'}\bigg] - \nu\frac{\varphi''\varphi - 2(\varphi')^2}{\varphi^3}w^2 - 2\varphi D_xu_\lambda \cdot D_xH + 2\gamma\varphi D_xu_\lambda \cdot D_xf \\
        &\leq C_{R,\varepsilon}\left(\Lambda_{q'}(\mu)^{q'} + 1\right) + \lambda\left(\|u_\lambda\|_\infty + (1 + C_0)\int_t^T e^{\lambda(t-s)}\Lambda_{q'}(\mu_\lambda)^{q'}ds\right)R^{-1}w \\
        &\qquad + C_0\gamma(2ew^\frac{1}{2} + R^{-1}w) - \varepsilon w^{1+\frac{q}{2}}
    \end{aligned}
    $$
    Choose $R \geq 2\left(\|u_\lambda\|_\infty + (1 + C_0)\underset{t}{\sup}\int_t^T e^{\lambda(t-s)}\Lambda_{q'}(\mu_\lambda)^{q'}ds + C_0\right)$. If $\gamma = 1$, we note that $C_0\gamma(2ew^\frac{1}{2} + R^{-1}w) \leq \varepsilon w^{1+\frac{q}{2}} + C_\varepsilon$ and $\left(\|u_\lambda\|_\infty + (1 + C_0)\int_t^T e^{\lambda(t-s)}\Lambda_{q'}(\mu_\lambda)^{q'}ds\right)R^{-1} \leq 1$. If $\gamma = \lambda$, we use that $\left(\|u_\lambda\|_\infty + (1 + C_0)\int_t^T e^{\lambda(t-s)}\Lambda_{q'}(\mu_\lambda)^{q'}ds + C_0\right)R^{-1} \leq \frac{1}{2}$ and $2eC_0w^\frac{1}{2} \leq \frac{1}{2}w + 16e^2C_0^2$. In either case, we get
    $$
    \begin{aligned}
        &-\partial_tw - \nu\Delta w + D_xw \cdot D_pH - 2\nu\frac{\varphi'}{\varphi}D_xw \cdot D_xu_\lambda + 2\nu\varphi|D_{xx}^2u_\lambda|^2 + \lambda w \\
        &\leq C_{R,\varepsilon}\left(\Lambda_{q'}(\mu)^{q'} + 1 + \gamma\right)
    \end{aligned}
    $$
    and so the maximum principle gives that
    $$
    \begin{aligned}
        \|w(t,\cdot)\|_\infty &\leq e^{\lambda(t-T)}\|w(T,\cdot)\|_\infty + C\int_t^T e^{\lambda(t-s)} \left(\Lambda_{q'}(\mu)^{q'} + 1 + \gamma\right)ds \\
        &\leq C\left((1+\gamma + \gamma^2)e^{\lambda(t-T)} + (1 + \gamma)\lambda^{-1}(1 - e^{\lambda(t-T)}) \right).
    \end{aligned}
    $$

    In the case where $\lambda < 1$, we choose $R \geq \lambda\left(\|u_\lambda\|_\infty + (1+C_0)\underset{t}{\sup}\int_t^T e^{\lambda(t-s)}\Lambda_{q'}(\mu_\lambda(s))^{q'} ds\right)$ and define
    $$w(t,x) = \lambda\varphi\left(\lambda u_\lambda(t,x) + (1+C_0)\lambda\int_t^T e^{\lambda(t-s)}\Lambda_{q'}(\mu_\lambda(s))^{q'} ds\right)|D_xu_\lambda(t,x)|^2.$$
    Since $\gamma \leq 1$, we have
    $$
    \begin{aligned}
        &-\partial_tw - \nu\Delta w + D_xw \cdot D_pH - 2\nu\lambda\frac{\varphi'}{\varphi}D_xw \cdot D_xu_\lambda + 2\nu\lambda\varphi|D_{xx}^2u_\lambda|^2 + 2\lambda w \\
        &= \lambda\frac{\varphi'}{\varphi}w\bigg[\gamma f - \lambda u_\lambda + D_xu_\lambda \cdot D_pH - H - \lambda(1+C_0)\int_t^T e^{\lambda(t-s)}\Lambda_{q'}(\mu_\lambda(s))^{q'} ds \\
        &\qquad\qquad + (1 + C_0)\Lambda_{q'}(\mu_\lambda)^{q'}\bigg] - \nu\frac{\varphi''\varphi - 2(\varphi')^2}{\varphi^3}w^2 - 2\lambda\varphi D_xu_\lambda \cdot D_xH + 2\gamma \lambda\varphi D_xu_\lambda \cdot D_xf \\
        &\leq \lambda\left(\varepsilon - \frac{1}{C}\right)\left(w^{1+\frac{q}{2}} + w\Lambda_{q'}(\mu)^{q'}\right) + \lambda C_{R,\varepsilon}\left(\Lambda_{q'}(\mu)^{q'} + 1\right),
    \end{aligned}
    $$
    and so the maximum principle yields
    $$
    \begin{aligned}
        \|w(t,\cdot)\|_\infty &\leq e^{\lambda(t-T)}\|w(T,\cdot)\|_\infty + C\lambda\int_t^T e^{\lambda(t-s)} \left(\Lambda_{q'}(\mu)^{q'} + 1 + \gamma\right)ds \\
        &\leq C\left((1+\gamma + \lambda\gamma)e^{\lambda(t-T)} + (1 + \gamma)(1 - e^{\lambda(t-T)}) \right) \\
        &\leq C.
    \end{aligned}
    $$
    Thus, in both cases, $|D_xu_\lambda| \leq \lambda^{-\frac{1}{2}}w^\frac{1}{2} \leq C\lambda^{-\frac{1}{2}}$, recalling that if $T = \infty$, the terminal condition is dropped and we have $\int_t^\infty e^{\lambda(t-s)} ds = \lambda^{-1}$.
\end{proof}

By Lemma \ref{Lem: General C1 Estimate}, \cite[Theorem 4.1]{kobeissi2022mean} yields estimates for $\Lambda_\infty(\mu_\lambda)$. From this, we can apply (5.6)-(5.9) in \cite{kobeissi2022mean} to obtain bounds for $H$, $D_pH$, and $D_xH$.

\begin{lemma} \label{Lem: Holder Estimate for m}
    Assume A\ref{A: A1}-\ref{A: A6} hold. Fix $\nu_0 \in (0,\infty)$, $\nu \in (0,\nu_0)$, $\lambda \geq \frac{1}{C_0}$, and $\gamma \in \{1,\lambda\}$. Now let $(u_\lambda,m_\lambda,\mu_\lambda)$ be a solution to \eqref{Eq: General Discounted MFG} in the sense of distributions, where $\Omega = \T^d$ and $T \in (0,\infty)$. Then there exists some $C > 0$ independent of $\nu,\lambda,T$ such that
    $$W_1(m_\lambda(t_1),m_\lambda(t_2)) \leq C|t_1-t_2|^{1/2}$$
    for all $t_1,t_2 \in [0,T]$ with $|t_1-t_2| \leq 1$.
\end{lemma}

\begin{proof}
    Note $X(t) \sim m_\lambda(t)$, where $dX(t) = -D_pH(t,X(t),D_xu_\lambda(t,X(t)),\mu_\lambda(t))dt + \sqrt{2\nu}dB(t)$, and so
    $$
    \begin{aligned}
        &W_1(m(t_1),m(t_2)) \\
        &= \sup_{Lip(\varphi) \leq 1} \int_{\T^d} \varphi(x) d(m(t_1) - m(t_2)) \\
        &= \sup_{Lip(\varphi) \leq 1} E|\varphi(X(t_1)) - \varphi(X(t_2))| \\
        &\leq E|X(t_1) - X(t_2)| \\
        &= E\left|-\int_{t_2}^{t_1} D_pH(s,X(s),D_xu_\lambda(s,X(s)),\mu_\lambda(s)) ds + \sqrt{2\nu}\int_{t_2}^{t_1} dB(s)\right| \\
        &\leq |t_1-t_2|^\frac{1}{2}\left[E\int_{t_2}^{t_1} |D_pH(s,X(s),D_xu_\lambda(s,X(s)),\mu_\lambda(s))|^2 ds\right]^\frac{1}{2} \\
        &\qquad+ \sqrt{2\nu d}\left[E\left|\int_{t_2}^{t_1} dB(s)\right|^2\right]^\frac{1}{2} \\
        &\leq |t_1-t_2|^\frac{1}{2}\left(\sqrt{2\nu d} + \left[\int_{t_2}^{t_1}\int_{\T^d} |D_pH(s,X(s),D_xu_\lambda(s,X(s)),\mu_\lambda(s))|^2 dm_\lambda(s)ds\right]^\frac{1}{2}\right) \\
        &\leq C|t_1-t_2|^\frac{1}{2}, \qquad \forall |t_1-t_2| \leq 1.
    \end{aligned}
    $$
\end{proof}

\begin{lemma} \label{Lem: W1 estimate}
    Suppose A\ref{A: A1}-\ref{A: A6} and B\ref{A: DpH is Lipschitz in alpha}-\ref{A: Regularity of DpH} hold, fix $0 \leq r \leq s < \infty$ and $p_1,p_2 \in C^{0,1}([r,s] \times \T^d; \R^d)$ with $\|p_1(t,\cdot)\|_{C^1(\T^d)} + \|p_2(t,\cdot)\|_{C^1(\T^d)} \leq C_0$, and let $(m_1,\mu_1)$, $(m_2,\mu_2)$ be solutions to
    $$
    \begin{cases}
        \partial_tm_i - \nu\Delta m_i - \operatorname{div}(m_iD_pH(t,x,p_i,\mu_i)) = 0, &(t,x) \in [r,s] \times \T^d \\
        \mu_i = (I,-D_pH(t,x,p_i,\mu_i)) \# m_i, &t \in [r,s]
    \end{cases}
    $$
    Then we have
    $$W_1(m_1(s),m_2(s)) \leq \left(W_1(m_1(r),m_2(r)) + C\int_r^s \|p_1(\cdot,t) - p_2(\cdot,t)\|_\infty dt\right)e^{C(s-r)}.$$
\end{lemma}

\begin{proof}
    Let $\xi_1,\xi_2$ be random variables in $\T^d$ with $m_i(r) = \mathcal{L}(\xi_i)$ and
    $$W_1(m_1(r),m_2(r)) = E|\xi_1 - \xi_2|.$$
    Now note that for $r \leq t \leq s$, $m_i(t) = \mathcal{L}(X_t^i)$ where
    $$
    \begin{cases}
        dX_t^i = -D_pH(t,X_t^i,p_i(t,X_t^i),\mu_i(t))dt + \sqrt{2\nu}dB_t \\
        X_r^i = \xi_i
    \end{cases}
    $$
    Thus, we have
    $$
    \begin{aligned}
        &W_1(m_1(s),m_2(s)) \leq E|X_s^1 - X_s^2| \\
        &\quad\leq E|X_r^1 - X_r^2| + \int_r^s E|D_pH(t,X_t^1,p_1(t,X_t^1),\mu_1(t)) - D_pH(t,X_t^2,p_2(t,X_t^2),\mu_2(t))| dt
    \end{aligned}
    $$
     Moreover, by B\ref{A: DpH is Lipschitz in alpha}-\ref{A: Regularity of DpH}, we find
    $$
    \begin{aligned}
        &(1 - L_1)|D_pH(t,X_t^1,p_1(t,X_t^1),\mu_1(t)) - D_pH(t,X_t^2,p_2(t,X_t^2),\mu_2(t))| \\
        &\qquad \leq C(|X_t^1 - X_t^2| + |p_1(t,X_t^1) - p_2(t,X_t^2)| + W_1(m_1(t),m_2(t))).
    \end{aligned}
    $$
    Therefore, we have
    $$
    \begin{aligned}
        E|X_s^1 - X_s^2| &\leq E|X_r^1 - X_r^2| \\
        &\qquad+ C\int_r^s E(|X_t^1 - X_t^2| + |p_1(t,X_t^1) - p_2(t,X_t^2)| + W_1(m_1(t),m_2(t))) dt \\
        &\leq E|X_r^1 - X_r^2| + C\int_r^s E(2|X_t^1 - X_t^2| + |p_1(t,X_t^1) - p_2(t,X_t^2)|) dt \\
        &\leq E|X_r^1 - X_r^2| + C\int_r^s \left[(2 + C_0)E|X_t^1 - X_t^2| + \|p_1(t,\cdot) - p_2(t,\cdot)\|_\infty\right] dt,
    \end{aligned}
    $$
    where the last inequality follows because
    $$
    \begin{aligned}
        |p_1(t,X_t^1) - p_2(t,X_t^2)| &\leq |p_1(t,X_t^1) - p_1(t,X_t^2)| + |p_1(t,X_t^2) - p_2(t,X_t^2)| \\
        &\leq C_0|X_t^1 - X_t^2| + \|p_1(t,\cdot) - p_2(t,\cdot)\|_\infty.
    \end{aligned}
    $$
    Thus, Gronwall's inequality yields
    $$
    \begin{aligned}
        &W_1(m_1(s),m_2(s)) \leq E|X_s^1 - X_s^2| \leq \left(E|X_r^1 - X_r^2| + C\int_r^s \|p_1(t,\cdot) - p_2(t,\cdot)\|_\infty dt\right)e^{\widetilde{C}(s-r)} \\
        &\qquad = \left(W_1(m_1(r),m_2(r)) + C\int_r^s \|p_1(\cdot,t) - p_2(\cdot,t)\|_\infty dt\right)e^{\widetilde{C}(s-r)}.
    \end{aligned}
    $$
\end{proof}

%-----%-----%-----%
%-----%-----%-----%
%-----%-----%-----%
%-----%-----%-----%
\subsection{Existence on Infinite Time Horizons}
\label{Sec: Existence}
%-----%-----%-----%
%-----%-----%-----%
%-----%-----%-----%
%-----%-----%-----%

\begin{lemma} \label{Lem: Existence}
    Assume A\ref{A: A1}-\ref{A: A6} hold. Then for $\nu,\lambda > 0$ and $\gamma \in \{1,\lambda\}$, there exists at least one solution $(u_\lambda,m_\lambda,\mu_\lambda)$ to \eqref{Eq: General Discounted MFG} in the sense of distributions, where $\Omega \in \{\T^d,\R^d\}, T \in (0,\infty)$. Furthermore, if B\ref{A: A1}-\ref{A: Regularity of H} hold, then there exists at least one classical solution to \eqref{Eq: General Discounted MFG}, where either $\Omega \in \{\T^d,\R^d\}, T \in (0,\infty)$ or $\Omega = \T^d,T \in (0,\infty]$ (the terminal condition is dropped when $T = \infty$).
\end{lemma}

\begin{proof}
    In the case $T < \infty$, the existence of solutions in the sense of distributions follows from \cite{kobeissi2022mean} because
    $$-\partial_t u_\lambda -\nu\Delta u_\lambda + H(t,x,D_xu_\lambda,\mu) + \lambda u_\lambda = \gamma f(t,x,m_\lambda(t))$$
    if and only if
    $$-\partial_t v_\lambda -\nu\Delta v_\lambda + e^{\lambda(T-t)}H(t,x,e^{\lambda(t-T)}D_xv_\lambda,\mu) = \gamma e^{\lambda(T-t)}f(t,x,m_\lambda(t))$$
    where $v_\lambda = e^{\lambda(T-t)}u_\lambda$.  Moreover, standard bootstrapping arguments (see \cite{kobeissi2022classical} or \cite[Section 2.5]{graber2025mean}) give $H(t,x,D_xu_\lambda,\mu_\lambda) \in C^{\alpha/2,\alpha}([0,T] \times \T^d)$ and hence $u_\lambda \in C^{1+\alpha/2,2+\alpha}([0,T] \times \T^d)$ for some $\alpha \in (0,1)$.

    Now suppose $T = \infty$ and fix $T_0 > 0$. For all $\tau > 0$, let $(u^\tau,m^\tau,\mu^\tau)$ be a solution to \eqref{Eq: General Discounted MFG} on $(0,\tau) \times \T^d$ with $g \equiv 0$. By Lemmas \ref{Lem: General C1 Estimate} and \ref{Lem: Holder Estimate for m}, we have $\|u^\tau\|_{C^{0,1}([0,\tau] \times \T^d)} + \|m^\tau\|_{C^{1/2}(0,T_0;\mathcal{P}(\T^d))} \leq C$ for all $\tau \geq T_0$, where $C$ may depend on $T_0$ but not on $\tau$. Applying \cite[Theorem IV.6.1]{ladyzhenskaia1968linear},
    $$
    \begin{aligned}
        \|u^\tau\|_{C^{1+\alpha/2,2+\alpha}([0,\tau] \times \T^d)} &\leq C(\nu)\left(\|H\|_{\alpha/2,\alpha} + \|f\|_{\alpha/2,\alpha} + \lambda\|u^\tau\|_{\alpha/2,\alpha}\right) \\
        &\leq C(\nu,\lambda)\left(\|H\|_{\alpha/2,\alpha} + 1\right).
    \end{aligned}
    $$
    Combining B\ref{A: Regularity of H} with the bounds on $D_pH$ and $D_xH$ (e.g., see Lemma \ref{Lem: General C1 Estimate} and \cite[Lemma 2.11]{kobeissi2022mean}), we have
    $$
    \begin{aligned}
        &|H(t_1,x,D_xu^\tau(t_1,x),\mu^\tau(t_1)) - H(t_2,y,D_xu^\tau(t_2,y),\mu^\tau(t_2))| \\
        &\leq C\left[\left([D_xu^\tau]_{\alpha/2,\alpha} + 1\right)(|t_1-t_2|^{\alpha/2} + |x-y|^\alpha) + W_1(\mu^\tau(t_1),\mu^\tau(t_2))\right], 
    \end{aligned}
    $$
    and applying B\ref{A: DpH is Lipschitz in alpha}-\ref{A: Regularity of DpH} yields
    $$
    \begin{aligned}
        W_1(\mu^\tau(t_1),\mu^\tau(t_2)) &= \sup_{Lip(\phi)} \bigg(\int_{\T^d} \phi(x,-D_pH(x,t_1,D_xu^\tau(t_1,x),\mu^\tau(t_1))) d(m^\tau(t_1) - m^\tau(t_2)) \\
        &\qquad\qquad+ \int_{\T^d} \big(\phi(x,-D_pH(x,t_1,D_xu^\tau(t_1,x),\mu^\tau(t_1))) \\
        &\qquad\qquad\qquad- \phi(x,-D_pH(x,t_2,D_xu^\tau(t_2,x),\mu^\tau(t_2)))\big) dm^\tau(t_2)\bigg) \\
        &\leq C\bigg[(\|D_{xx}^2u^\tau\|_\infty + 1)W_1(m^\tau(t_1),m^\tau(t_2)) + ([D_xu]_{\alpha/2,\alpha} + 1)|t_1 - t_2|^{\alpha/2}\bigg].
    \end{aligned}
    $$
    Moreover, Lemma \ref{Lem: Holder Estimate for m} yields $W_1(m^\tau(t_1),m^\tau(t_2)) \leq C|t_1-t_2|$ for $|t_1-t_2| \leq 1$, and the boundedness of $H$ gives us that this holds for all $t_1,t_2 \in (0,\tau)$. By interpolation on H\"{o}lder spaces, we get
    $$\|D_{xx}^2u^\tau\|_\infty + [D_xu^\tau]_{\alpha/2,\alpha} \leq \delta \|u^\tau\|_{C^{1+\alpha/2,2+\alpha}([0,\tau] \times \T^d)} + C_\delta\|u^\tau\|_\infty$$
    for $\delta > 0$.
    Finally, choosing $\delta > 0$ sufficiently small, we get a bound for $\|u^\tau\|_{C^{1+\alpha/2,2+\alpha}([0,T_0] \times \T^d)}$ that does not depend on $\tau$.
    
    By the Arzela-Ascoli Theorem and a diagonalization argument, we get that $(u^\tau,m^\tau)$ has a subsequence (which we will not rename) that converges to some $(u,m)$ in $C^{1,2}([0,T_0] \times \T^d) \times C^0(0,T_0;\mathcal{P}(\T^d))$ for every $T_0 > 0$. By \cite[Lemma 4.4]{kobeissi2022mean}, it follows that $\mu^\tau \to \mu$ in $C^0(0,T_0;\mathcal{P}(\T^d))$ for every $T_0 > 0$, where $\mu = (I,-D_pH(t,\cdot,D_xu,\mu)) \# m$. Passing to the limit, we get that $(u,m,\mu)$ satisfies \eqref{Eq: General Discounted MFG}, thus completing the proof.
\end{proof}

%-----%-----%-----%
%-----%-----%-----%
%-----%-----%-----%
%-----%-----%-----%
\subsection{Asymptotic Behavior of a Discounted MFGC}
\label{Sec: ABM}
%-----%-----%-----%
%-----%-----%-----%
%-----%-----%-----%
%-----%-----%-----%

\begin{remark}
    For $\lambda > 0$, Lemma \ref{Lem: Existence} gives that there exists a solution $(u_\lambda,m_\lambda,\mu_\lambda)$ to \eqref{Eq: MFG 1} in the sense of distributions.
    Furthermore, by Lemma \ref{Lem: General C1 Estimate}, we have $\|u_\lambda\|_\infty + \|D_xu_\lambda\|_\infty \leq C(e^{\frac{\lambda}{2}(t-T)} + \lambda^{-\frac{1}{2}})$ for $\lambda \geq 1$, which gives us that $u_\lambda(t,\cdot)$ converges to $0$ in $C^1(\T^d)$ locally uniformly in $[0,T)$ as $\lambda \to \infty$.
\end{remark}

\begin{theorem} \label{Thm: MFG to ABM}
    Assume A\ref{A: A1}-\ref{A: A6} hold. Fix $T \in (0,\infty)$ and let $(u_\lambda,m_\lambda,\mu_\lambda)_\lambda$ be a sequence of solutions to \eqref{Eq: MFG 1} with $\nu = \nu_\lambda$ where $(\nu_\lambda)_\lambda$ is a sequence of positive numbers with $\nu_\lambda \downarrow \nu_\infty \geq 0$ as $\lambda \to \infty$. Then as $\lambda \to \infty$, we get that $u_\lambda \to 0$ in $C^{0,1}([0,T-\varepsilon] \times \T^d)$ for every $\varepsilon > 0$, $(m_\lambda)_\lambda$ is bounded in $C^{1/2}([0,T]; \mathcal{P}(\T^d))$ for some $\alpha \in (0,1)$, and any limit point $m_*$ of $(m_\lambda)_\lambda$ in $C^0(0,T;\mathcal{P}(\T^d))$ as $\lambda \to \infty$ is a solution to the model \eqref{Eq: ABM}
    \begin{equation*}
        \begin{cases}
            \partial_t m_* - \nu_\infty\Delta m_* - \operatorname{div} (m_*D_pH(t,x,0,\mu_*)) = 0, &(t,x) \in (0,T) \times \T^d \\
            m_*(0,x) = m_0(x), &x \in \T^d
        \end{cases}
    \end{equation*}
    in the sense of distributions, where
    $$\mu_* = (I,-D_pH(t,\cdot,0,\mu_*)) \# m_*, \qquad \forall t \in (0,T).$$
\end{theorem}

\begin{proof}
    The convergence of $(u_\lambda)_\lambda$ follows from Lemma \ref{Lem: General C1 Estimate}. By Lemma \ref{Lem: Holder Estimate for m}, we have that $(m_\lambda)_\lambda$ is a uniformly bounded and equicontinuous sequence in $C^0(0,T;\mathcal{P}(\T^d))$.

    Now suppose $m_\lambda \to m_*$ in $C^0(0,T;\mathcal{P}(\T^d))$ (up to a subsequence which we will not rename). By \cite[Lemma 4.4]{kobeissi2022mean}, we get that $\mu_\lambda(t) \to \mu_*(t)$ for all $t < T$. Further, for $\varphi \in C^\infty_c([0,T) \times \T^d)$, we have
    $$\int_0^T\int_{\T^d} (-\partial_t\varphi -\nu_\lambda\Delta\varphi + D_pH(t,x,D_xu_\lambda,\mu_\lambda) \cdot D_x\varphi) dm_\lambda(t) dt = \int_{\T^d} \varphi(0,\cdot)dm_0.$$
    Passing to the limit as $\lambda \to \infty$, we get
    $$\int_0^T\int_{\T^d} (-\partial_t\varphi -\nu_\infty\Delta\varphi + D_pH(t,x,0,\mu_*) \cdot D_x\varphi) dm_*(t) dt = \int_{\T^d} \varphi(0,\cdot)dm_0$$
    since $D_pH(t,x,D_xu_\lambda,\mu_\lambda) \to D_pH(t,x,0,\mu_*)$ locally uniformly and is uniformly bounded (see Lemma \ref{Lem: General C1 Estimate}).
\end{proof}

\begin{remark}
    The case where $T = \infty$ is very similar, except that we require assumptions B\ref{A: DpH is Lipschitz in alpha}-\ref{A: Regularity of H} to guarantee the existence of solutions for all $\lambda$. In this case, we get that $(m_\lambda)_\lambda$ is bounded in $C^{\alpha/2,\alpha}([0,T_0] \times \T^d)$ for all $T_0$ and use a diagonalization argument to extract a subsequence that converges in $C^0(0,T_0;\mathcal{P}(\T^d))$ for all $T_0 > 0$.
\end{remark}

\begin{remark}
    If we assume that A\ref{A: A1}-\ref{A: A6} and B\ref{A: DpH is Lipschitz in alpha}-\ref{A: Regularity of DpH} hold (and B\ref{A: Regularity of H} if $T = \infty$), and $D_pH$ does not depend on $t$, then the vector field $G(x,m) \coloneqq -D_pH(x,0,\mu_m)$ (where $\mu_m = (I,-D_pH(x,0,\mu_m) \# m$) satisfies
    $$|G(x,m) - G(y,m)| \leq C|x-y|$$
    by B\ref{A: Regularity of DpH}, and
    $$|G(x,m_1) - G(x,m_2)| \leq CW_1(\mu_{m_1},\mu_{m_2}) \leq CW_1(m_1,m_2)$$
    since
    $$
    \begin{aligned}
        W_1(\mu_{m_1},\mu_{m_2}) &= \sup_{Lip(\varphi)} \Bigg[\int_\Omega (\varphi(x,-D_pH(x,0,\mu_{m_1})) - \varphi(x,-D_pH(x,0,\mu_{m_2}))) dm_1 \\
        &\qquad + \int_{\T^d} \varphi(x,-D_pH(x,0,\mu_{m_2})) d(m_1 - m_2)\Bigg] \\
        &\leq CW_1(m_1,m_2)
    \end{aligned}
    $$
    (where the last inequality follows from B\ref{A: DpH is Lipschitz in alpha}-\ref{A: Regularity of DpH}). Thus, under these assumptions, one can see that \eqref{Eq: ABM} has a unique solution by standard Cauchy-Lipschitz theory (e.g., see \cite{bogachev2022fokker,piccoli2013transport}), and hence the whole sequence $(u_\lambda,m_\lambda,\mu_\lambda)_\lambda$ converges to a unique limit as $\lambda \to \infty$.
\end{remark}

%-----%-----%-----%
%-----%-----%-----%
%-----%-----%-----%
%-----%-----%-----%
\subsection{Asymptotic Behavior of the Rescaled MFGC}
\label{Sec: Aggregation}
%-----%-----%-----%
%-----%-----%-----%
%-----%-----%-----%
%-----%-----%-----%

\begin{remark}
    For $\lambda > 0$, $(u_\lambda,m_\lambda,\mu_\lambda)$ solves \eqref{Eq: Rescaled MFG} if and only if $(v_\lambda,m_\lambda,\mu_\lambda)$ solves
    \begin{equation} \label{Eq: COV 1}
    \begin{cases}
        -\partial_t v_\lambda -\nu_\lambda\Delta v_\lambda + H(t,x,D_xv_\lambda,\mu_\lambda) + \lambda v_\lambda = \lambda f(t,x,m(t)), &(t,x) \in (0,\infty) \times \T^d \\
        \partial_t m_\lambda - \nu_\lambda\Delta m_\lambda - \operatorname{div} (m_\lambda D_pH(t,x,D_x v_\lambda,\mu_\lambda)) = 0, &(t,x) \in (0,\infty) \times \T^d \\
        \mu_\lambda = (I,-D_pH(t,\cdot,D_xv_\lambda(t,\cdot),\mu_\lambda)) \# m_\lambda, &t \in (0,\infty) \\
        m_\lambda(0,x) = m_0(x), &x \in \T^d
    \end{cases}
\end{equation}
where $v_\lambda = \lambda u_\lambda$. By Lemma \ref{Lem: Existence}, there exist classical solutions to \eqref{Eq: COV 1}. Furthermore, Lemma \ref{Lem: General C1 Estimate} gives us that $\|v_\lambda\|_\infty + \|D_xv_\lambda\|_\infty \leq C$. In particular, $\|u_\lambda(t,\cdot)\|_\infty + \|D_xu_\lambda(t,\cdot)\|_\infty \leq C\lambda^{-1}$.
\end{remark}

\begin{lemma} \label{Lem: Semiconcavity of u}
    Assume A\ref{A: A1}-A\ref{A: A6}, B\ref{A: DpH is Lipschitz in alpha}-B\ref{A: Regularity of H}, and C\ref{A: f is semiconcave}-C\ref{A: Second Derivatives of H} hold and that $\lambda \geq 1$. If $(u_\lambda,m_\lambda,\mu_\lambda)$ is a solution of \eqref{Eq: Rescaled MFG}, then $D_{xx}^2u_\lambda \leq C\lambda^{-1}$.
\end{lemma}

\begin{proof}
    We will again assume $u_\lambda$ is smooth (if not, standard approximation arguments can be used). Fixing $\xi \in \R^d$ with $|\xi| \leq 1$, we let $w = D_xu_\lambda \cdot \xi$ and $z = D_{xx}^2u_\lambda \xi \cdot \xi$. Then we have
    $$
    \begin{aligned}
        -\partial_tz - \nu_\lambda\Delta z + \lambda z + D_pH \cdot D_xz + 2D_{p\xi}^2 H \cdot D_xw + \lambda D_{pp}^2H D_xw \cdot D_xw + \lambda^{-1}D_{\xi\xi}^2H = D_{\xi\xi}^2f
    \end{aligned}
    $$
    By C\ref{A: f is semiconcave}, we get $D_{\xi\xi}^2f \leq C$. Combining B\ref{A: Regularity of DpH}, C\ref{A: Second Derivatives of H}, and the uniform bound on $\lambda|D_xu_\lambda|$, we get that
    $$\lambda D_{pp}^2H(t,x,\lambda D_xu_\lambda,\mu_\lambda) D_xw \cdot D_xw \geq \frac{\lambda}{\widetilde{C}}|D_xw|^2$$
    and
    $$2D_{p\xi}^2H(t,x,\lambda D_xu_\lambda,\mu_\lambda) \cdot D_xw + \lambda^{-1}D_{\xi\xi}^2H(t,x,\lambda D_xu_\lambda,\mu_\lambda) \geq -\frac{\lambda}{2\widetilde{C}}|D_xw|^2 - C\lambda^{-1}$$
    for $\lambda \geq 1$. Thus, we conclude that
    $$-\partial_tz - \nu_\lambda\Delta z + \lambda z + D_pH \cdot D_xz \leq C.$$
    By the maximum principle, we get
    $$z \leq C\int_t^\infty e^{\lambda(t-s)} ds = C\lambda^{-1}.$$
    And the preceding argument can be flipped to yield the full estimate $\abs{z}\leq C\lambda^{-1}$. Taking the supremum over all $\xi \in \R^d$ with $|\xi| \leq 1$ concludes the proof.
\end{proof}

\begin{lemma} \label{Lem: Moment Estimates for m}
    Assume A\ref{A: A1}-A\ref{A: A6}, B\ref{A: DpH is Lipschitz in alpha}-B\ref{A: Regularity of H}, and C\ref{A: f is semiconcave}-C\ref{A: Second Derivatives of H} hold and that $\lambda \geq 1$. If $(u_\lambda,m_\lambda,\mu_\lambda)$ is a solution of \eqref{Eq: Rescaled MFG}, then for any $T > 0$, $(m_\lambda)_\lambda$ is relatively compact in $C^0(0,T;\mathcal{P}(\T^d))$ and bounded in $L^\infty([0,T] \times \T^d)$.
\end{lemma}

\begin{proof}
    Compactness follows from Lemma \ref{Lem: Holder Estimate for m}.
    For the $L^\infty$ estimate, we will again assume $u_\lambda$ and $m_\lambda$ are smooth. Then
    $$\partial_tm_\lambda - \nu_\lambda\Delta m_\lambda - m_\lambda \operatorname{Tr}\left(\lambda D_{pp}^2H D_{xx}^2u_\lambda + D_{px}^2H\right) - D_pH \cdot D_xm_\lambda = 0.$$
    By B\ref{A: Regularity of DpH} and the bounds on $\lambda D_xu_\lambda$, we get that $\operatorname{Tr}(D_{px}^2H(t,x,\lambda D_xu_\lambda,\mu_\lambda)) \leq C$. By the same reasoning, we have $|D_{pp}^2H(t,x,\lambda D_xu_\lambda,\mu_\lambda)| \leq C$. Since $H$ is convex and $\lambda D_{xx}^2u_\lambda \leq C$, this implies $\lambda\operatorname{Tr}(D_{pp}^2H D_{xx}^2u_\lambda) \leq C$.
    In particular,
    $$\partial_tm_\lambda - \nu_\lambda\Delta m_\lambda - Cm_\lambda - D_pH(t,x,\lambda D_xu_\lambda,\mu_\lambda) \cdot D_xm_\lambda \leq 0.$$
    Thus, the maximum principle gives us that
    $$\|m_\lambda(t,\cdot)\|_\infty \leq e^{Ct}\|m_0\|_\infty \leq C_0e^{CT}
    < \infty.$$
\end{proof}

\begin{theorem} \label{Thm: Convergence of Rescaled MFGC}
    Assume that A\ref{A: A1}-A\ref{A: A6}, B\ref{A: DpH is Lipschitz in alpha}-B\ref{A: Regularity of H}, and C\ref{A: f is semiconcave}-C\ref{A: Second Derivatives of H} hold and that $(\nu_\lambda)_{\lambda \geq 1}$ is a non-increasing sequence of positive numbers with $\nu_\lambda \downarrow \nu_\infty \geq 0$ as $\lambda \to \infty$. If $(u_\lambda,m_\lambda,\mu_\lambda)$ solves \eqref{Eq: Rescaled MFG} for each $\lambda$, then $$\underset{t}{\sup}\|u_\lambda(t,\cdot)\|_{C^1(\T^d)} \leq C\lambda^{-1},$$ the sequence
    $(m_\lambda)_\lambda$ is relatively compact in $C^0(0,T;\mathcal{P}(\T^d))$ and bounded in $L^\infty([0,T] \times \T^d)$ for every $T > 0$, and any limit point $m_*$ of $(m_\lambda)_\lambda$ as $\lambda \to \infty$ satisfies \eqref{Eq: Aggregation}:
    \begin{equation*} 
        \begin{cases}
            \partial_t m_* - \nu_\infty\Delta m_* - \operatorname{div} (m_*D_pH(t,x,D_xf(t,x,m_*(t)),\mu_*)) = 0, &(t,x) \in (0,\infty) \times \T^d \\
            m_*(0,x) = m_0(x), &x \in \T^d
        \end{cases}
    \end{equation*}
    where $\mu_* = (I,-D_pH(t,\cdot,D_xf(t,\cdot,m_*(t)),\mu_*)) \# m_*$.
\end{theorem}

\begin{proof}
    Given Lemmas \ref{Lem: General C1 Estimate} and \ref{Lem: Moment Estimates for m}, what remains is to show the convergence result. Now suppose that, up to a subsequence, $m_{\lambda} \to m_*$ in $C^0(0,T;\mathcal{P}(\T^d))$ for all $T > 0$. By A\ref{A: A6}, we have $f(\cdot,\cdot,m_{\lambda}) \to f(\cdot,\cdot,m_*)$ in $C^{0,1}([0,T] \times \T^d)$ for all $T > 0$. Since $w_{\lambda} \coloneqq \lambda u_{\lambda}$ satisfies
    $$-\lambda^{-1}\partial_tw_{\lambda} - \lambda^{-1}\nu_{\lambda} \Delta w_{\lambda} + w_{\lambda} + \lambda^{-1}H(t,x,D_xw_{\lambda},\mu_{\lambda}) = f(t,x,m_{\lambda}),$$ 
    and $\abs{H(t,x,Dw_{\lambda},\mu_{\lambda})}\leq C$ (from the gradient estimates and \cite[Lemma 4.1]{kobeissi2022mean}) independent of $\lambda$,
    \begin{equation}
        -\lambda^{-1}\partial_tw_{\lambda} - \lambda^{-1}\nu_{\lambda} \Delta w_{\lambda} + w_{\lambda} - C\lambda^{-1} \leq f(t,x,m_{\lambda}).
    \end{equation}
    Taking $\underset{\lambda \to \infty}{\liminf}$ of both sides yields that $w_* = \underset{\lambda \to \infty}{\limsup} \, w_{\lambda}$ is a sub-solution to 
    \begin{equation}\label{eq: w}
       w = f(t,x,m_*). 
    \end{equation}
    and similarly $w^* = \underset{\lambda \to \infty}{\liminf} \, w_{\lambda}$ is a super-solution to \eqref{eq: w}. In particular, $w^* \leq f(t,x,m_*) \leq w^*$, which implies that $w_\lambda \to f(t,x,m_*)$ locally uniformly.
    
    By Lemmas \ref{Lem: General C1 Estimate} and \ref{Lem: Semiconcavity of u}, we get that $(\lambda u_{\lambda})_{\lambda}$ is uniformly bounded and uniformly semiconcave. Thus, by an argument nearly identical to the proof of \cite[Theorem 3.3.3]{cannarsa2004semiconcave}, we get that $(\lambda D_xu_{\lambda})_{\lambda}$ converges to $D_xf(t,x,m_*)$ a.e.
    By tightness, for each $t$, we can extract a further subsequence $(\mu_{\lambda_k})_{\lambda_k}$ converging to some $\mu \in \mathcal{P}(\T^d \times \R^d)$. Note that for a continuous function $\phi$,
    $$
    \begin{aligned}
        &\left|\int_{\T^d \times \R^d} \phi(x,\alpha) d(\mu_{\lambda_k}(t) - (I,-D_pH(t,x,D_xf(t,x,m_*(t)),\mu)) \# m_*(t))\right| \\
        &\leq \int_{\T^d} \left|\phi(x,-D_pH(t,x,\lambda_kD_xu_{\lambda_k}(t,x),\mu_{\lambda_k}(t))) - \phi(x,-D_pH(t,x,D_xf(t,x,m_*(t)),\mu))\right|dm_{\lambda_k} \\
        &\qquad + \left|\int_{\T^d} \phi(x,-D_pH(t,x,D_xf(t,x,m_*(t)),\mu))d(m_{\lambda_k} - m_*(t))\right| \\
        &\to 0
    \end{aligned}
    $$
    where the last step follows from the dominated convergence theorem and the fact that $m_{\lambda_k}(t) \to m_*(t)$ in $\mathcal{P}(\T^d)$. In particular, $\mu = (I,-D_pH(t,x,D_xf(t,x,m_*(t)),\mu)) \# m_*(t)$. Thus, by uniqueness of the fixed-point (see \cite[Lemma 4.3]{kobeissi2022mean}), we conclude that $\mu_\lambda(t) \to \mu_*(t)$ for all $t$.
    
    Finally, as in the proof of Theorem \ref{Thm: MFG to ABM}, we note that for $\varphi \in C^\infty_c([0,\infty) \times \T^d)$, we have
    $$\int_0^\infty\int_{\T^d} (-\partial_t\varphi -\nu_{\lambda}\Delta\varphi + D_pH(t,x,\lambda D_xu_{\lambda},\mu_{\lambda}) \cdot D_x\varphi) dm_{\lambda}(t) dt = \int_{\T^d} \varphi(0,\cdot)dm_0,$$
    and letting $\lambda \to \infty$ gives that
    $$\int_0^\infty\int_{\T^d} (-\partial_t\varphi -\nu_\infty\Delta\varphi + D_pH(t,x,D_xf(t,x,m_*),\mu_*) \cdot D_x\varphi) dm_*(t) dt = \int_{\T^d} \varphi(0,\cdot)dm_0.$$
    by the dominated convergence theorem.
\end{proof}

\begin{remark}
    Under the additional assumption that $D_pH$ does not depend on $t$ and
    \begin{equation}
        |D_xf(x,m_1) - D_xf(y,m_2)| \leq C_0(|x-y| + W_1(m_1,m_2))
    \end{equation}
    for $x,y \in \T^d$ and $m_1,m_2 \in \mathcal{P}(\T^d)$, letting $\mu_m = (I,-D_pH(x,D_xf(x,m),\mu_m) \# m$, we find that the vector field $G(x,m) \coloneqq -D_pH(x,D_xf(x,m),\mu_m)$ satisfies
    $$|G(x,m) - G(y,m)| \leq C(|x - y| + |D_xf(x,m) - D_xf(y,m)|) \leq C|x-y|$$
    by B\ref{A: Regularity of DpH}, and
    $$|G(x,m_1) - G(x,m_2)| \leq C(|D_xf(x,m_1) - D_xf(x,m_2)| + W_1(\mu_{m_1},\mu_{m_2})) \leq CW_1(m_1,m_2)$$
    since
    $$
    \begin{aligned}
        &W_1(\mu_{m_1},\mu_{m_2}) \\
        &= \sup_{Lip(\varphi)} \Bigg[\int_\Omega (\varphi(x,-D_pH(x,D_xf(x,m_1),\mu_{m_1})) - \varphi(x,-D_pH(x,D_xf(x,m_2),\mu_{m_2}))) dm_1 \\
        &\qquad + \int_{\T^d} \big(\varphi(x,-D_pH(x,D_xf(x,m_2),\mu_{m_2})) - \varphi(0,0)\big) d(m_1 - m_2)\Bigg] \\
        &\leq CW_1(m_1,m_2)
    \end{aligned}
    $$
    as before. Thus, under this additional assumption, we find that \eqref{Eq: Aggregation} has a unique solution and hence the whole sequence $(u_\lambda,m_\lambda,\mu_\lambda)_\lambda$ converges to a unique limit as $\lambda \to \infty$.
\end{remark}

\begin{remark}
    Suppose $\nu_\lambda = \nu_0 > 0$ for all $\lambda$, A\ref{A: A1}-\ref{A: A6} and B\ref{A: DpH is Lipschitz in alpha}-\ref{A: Regularity of H} hold, and $f: \mathcal{P}(\T^d) \to C^{2+\beta}(\T^d)$ with $\underset{m}{\sup} \|f\|_{C^2(\T^d)} \leq C_0$ and
    $$|f(x,m_1) - f(x,m_2)| + |D_xf(x,m_1) - D_xf(x,m_2)| \leq CW_1(m_1,m_2)$$
    for all $m_1,m_2 \in \mathcal{P}(\T^d)$. Now suppose $H$ does not depend on $t$ and is of class $C^{2+\beta}$ with respect to the variables $(x,p)$ with derivatives continuous on $\T^d \times \R^d \times \mathcal{P}(\T^d \times \R^d)$,
    %$$\zeta_1(|p|,\Lambda_\infty(\mu))^{-1}I \leq D_{pp}^2H(x,p,\mu) \leq \zeta_1(|p|,\Lambda_\infty(\mu))I,$$
    $$\frac{1}{C_R}I \leq D_{pp}^2H(x,p,\mu) \leq C_RI,$$
    and
    %$$|D_{px}^2H(x,p,\mu)| + |D_{xx}^2H(x,p,\mu)| \leq \zeta_1(|p|,\Lambda_\infty(\mu)),$$
    $$|D_{px}^2H(x,p,\mu)| + |D_{xx}^2H(x,p,\mu)| \leq C_R,$$
    for $|p| \leq R$ and $\mu \in \mathcal{P}_{\infty,R}(\T^d \times \R^d)$. Then by a nearly identical proof to the one found in \cite{cirant2025uniqueness}, we find that for $\lambda$ sufficiently large, we have
    $$\|D_{xx}^2u_\lambda\|_\infty \leq C\lambda^{-1} \qquad \text{and} \qquad \|\lambda D_xu_\lambda - D_xf(\cdot,m_*)\|_\infty \leq C\lambda^{-\frac{1}{2}}.$$
    Applying Lemma \ref{Lem: W1 estimate}, we find
    $$W_1(m_\lambda(t),m_*(t)) \leq Cte^{Ct}\lambda^{-\frac{1}{2}}$$
    and hence
    $$
    \begin{aligned}
        W_1(\mu_\lambda(t),\mu_*(t)) &\leq C(W_1(m_\lambda(t),m_*(t)) + \|\lambda D_xu_\lambda(t,\cdot) - D_xf(\cdot,m_*(t))\|_\infty) \\
        &\leq C(1 + te^{Ct})\lambda^{-\frac{1}{2}}.
    \end{aligned}
    $$
\end{remark}

%-----%-----%-----%
%-----%-----%-----%
%-----%-----%-----%
%-----%-----%-----%
\section{Stationary Models}
\label{Sec: Stationary}
%-----%-----%-----%
%-----%-----%-----%
%-----%-----%-----%
%-----%-----%-----%

In this section, we study the well-posedness of \eqref{Eq: Stationary Discounted MFGC}, \eqref{Eq: Ergodic MFGC}, and \eqref{Eq: Stationary First-Order MFGC}, as well as the convergence of solutions to \eqref{Eq: Stationary Discounted MFGC} as $\lambda \to 0$ or $\nu \to 0$. We begin by recalling standard results for the Hamilton-Jacobi and Fokker-Planck equations:

\begin{lemma} \label{Lem: Stationary HJ}
    Assume A\ref{A: A1}-\ref{A: A6} and D\ref{A: Bounds for H} hold. Then for each $\mu \in \mathcal{P}(\T^d \times \R^d)$, $\lambda \in (0,1)$, and $\nu \in (0,\frac{1}{d})$, there is a unique $u \in C^2(\T^d)$ satisfying
    $$-\nu\Delta u + H(x,D_xu,\mu) + \lambda u = f(x,\mu_1)$$
    (here $\mu_1$ is the first marginal of $\mu$). Furthermore,
    $$\lambda\|u\|_{L^\infty(\T^d)} \leq \| H(x,0,\mu) - f(x,\mu_1)\|_{L^\infty(\T^d)},$$
    and if C\ref{A: f is semiconcave} and D\ref{A: DpxH = 0 or H is uniformly convex} hold, then $D_{xx}^2u \leq C$, where $C = C(\|H\|_{C^1(\T^d)},\|f\|_{W^1_\infty(\T^d)}, \operatorname{SC}(f))$.
\end{lemma}

\begin{lemma} \label{Lem: Stationary FP}
    For every $b \in L^\infty(\T^d;\R^d)$ and $\lambda, \nu \in (0,1)$, there exists a unique weak solution $m \in H^1(\T^d)$ to
    $$-\nu\Delta m - \operatorname{div} (mb) + \lambda m = \lambda$$
    in the sense that
    $$\int_{\T^d} (\lambda m\varphi  + mb \cdot D_x\varphi + \nu D_xm \cdot D_x\varphi) dx = \lambda\int_{\T^d} \varphi dx \qquad \forall\varphi \in H^1(\T^d).$$
    Furthermore, $m \geq 0$, $\int_{\T^d} m dx = 1$, the set $\{x \in \T^d : m(x) = 0\}$ has measure zero, and for $d \geq 2$ there exists $C = C(d,\nu,\|b\|_\infty) > 0$ such that $\|m\|_{L^\infty(\T^d)} + \|m\|_{H^1(\T^d)} \leq C$.
\end{lemma}

To prove the existence of solutions, we will apply the Leray-Schauder fixed point theorem to a map $T$ which we will define as follows (see \cite{kobeissi2022mean}).

\begin{definition} \label{Def: Fixed-Point Map}
    First, define the space $X \coloneqq C^2(\T^d) \times L^2(\T^d)$. Now, for $\theta \in [0,1]$, define
    $$L^\theta(x,\alpha,\widetilde{\mu}) =
    \begin{cases}
        \theta L\left(x,\frac{\alpha}{\theta},(I \otimes \theta^{-1}I) \# \widetilde{\mu}\right), &0 < \theta \leq 1 \\
        0, &\theta = 0, \alpha = 0 \\
        \infty, &\theta = 0, \alpha \neq 0
    \end{cases}$$
    for which the associated Hamiltonian is
    $$H^\theta(x,p,\widetilde{\mu}) =
    \begin{cases}
        \theta H(x,p,(I \otimes \theta^{-1}I) \# \widetilde{\mu}), &0 < \theta \leq 1 \\
        0, &\theta = 0
    \end{cases}$$
    For $(\widetilde{u},\widetilde{m}) \in X$ and $\theta \in [0,1]$, define $\mu = (I,-D_pH^\theta(x,D_x\widetilde{u},\mu)) \# \eta(\widetilde{m})$ where
    $$\eta(\phi)(x) = \frac{\phi(x)_+ - \int_{\T^d}\phi(y)_+ dy}{\max\{1,\int_{\T^d}\phi(y)_+ dy\}} + 1,$$
    and define $(u,m) \in X$ by
    \begin{equation}
        \begin{cases} \label{Eq: MFG for Fixed-Point Map}
            -\nu\Delta u + H^\theta(x,D_xu,\mu) + \lambda u = \theta f(x,\eta(\widetilde{m})), &x \in \T^d \\
            -\nu\Delta m - \operatorname{div} (mD_pH^\theta(x,D_xu,\mu)) + \lambda m = \lambda, &x \in \T^d
        \end{cases}
    \end{equation}
    We define $T: X \times [0,1] \to X$ by $T(\widetilde{u},\widetilde{m},\theta) = (u,m)$.
\end{definition}

\begin{remark}
    The map $\eta$ found in Definition \ref{Def: Fixed-Point Map} was introduced in \cite{bonnans2019schauder}. Notice that $\eta$ is continuous. Moreover, for any $\phi \in L^1(\T^d)$,
    $$\eta(\phi)(x) = \frac{\phi(x)_+ - \int_{\T^d}\phi(y)_+ dy}{\max\{1,\int_{\T^d}\phi(y)_+ dy\}} + 1 \geq \frac{\phi(x)_+}{\max\{1,\int_{\T^d}\phi(y)_+ dy\}} \geq 0,$$
    $$\int_{\T^d} \eta(\phi)(x) dx = 1,$$
    and if $\phi \geq 0$ with $\int_{\T^d} \phi dx = 1$, then $\eta(\phi) = \phi$.
\end{remark}

\begin{lemma} \label{Lem: Estimates for Stationary Case}
    Assume A\ref{A: A1}-\ref{A: A6} and D\ref{A: Bounds for H} hold, $\lambda \in (0,1)$, and $\nu \in (0,\frac{1}{d})$. Then there exists some $C > 0$ (independent of $\nu,\theta$) such that if $(u,m) = T(u,m,\theta)$, then
    $$\Lambda_{q'}(\mu) + \|D_xu\|_\infty \leq C.$$
    Moreover, there exists some $\alpha \in (0,1)$ and $C(\nu) > 0$ independent of $\theta$ such that
    $$\|u - u(0)\|_{C^{2+\alpha(\T^d)}} \leq C(\nu).$$
\end{lemma}

\begin{proof}
    Note that
    \begin{equation} \label{Eq: Optimality of -DpH}
    \begin{aligned}
        \int_{\T^d \times \R^d} (L^\theta(x,\alpha,\mu) + \theta f(x,m)) d\mu &= \int_{\T^d} (D_xu \cdot D_pH^\theta(x,D_xu,\mu) - \nu\Delta u + \lambda u) dm \\
        &= \lambda\int_{\T^d} u dx \\
        &= \lambda\int_{\T^d} \inf_\alpha E\int_0^\infty (L^\theta(X_t,\alpha(X_t),\mu) + \theta f(X_t,m)) e^{-\lambda t} dt dx
    \end{aligned}
    \end{equation}
    Consider
    $$d\widetilde{X}_t = \sqrt{2\nu}dB_t, \qquad \widetilde{X}_0(x) = x$$
    and define $\widetilde{\mu}_t = \mathcal{L}(\widetilde{X}_t) \otimes \delta_0$.
    By A\ref{A: A3}-\ref{A: A4}, and noting A\ref{A: A3} holds for $L^\theta$ by a change of variables, we get
    $$
    \begin{aligned}
        \Lambda_{q'}(\mu)^{q'} &= \lambda\int_{\T^d}\int_0^\infty\int_{\T^d \times \R^d} |\alpha|^{q'}e^{-\lambda t} d\mu dtdx \\
        &\leq \lambda C_0\theta^{q'-1}\int_{\T^d}\int_0^\infty\int_{\T^d \times \R^d} L^\theta(x,\alpha,\widetilde{\mu}_t)e^{-\lambda t} d\mu dtdx + C_0^2\theta^{q'} \\
        &\leq \lambda C_0\theta^{q'-1}\int_{\T^d}\int_0^\infty\left[\int_{\T^d \times \R^d} L^\theta(x,\alpha,\mu) d(\mu - \widetilde{\mu}_t) + \int_{\T^d \times \R^d} L^\theta(x,\alpha,\widetilde{\mu}_t) d\widetilde{\mu}_t\right]e^{-\lambda t} dtdx \\
        &\qquad + C_0^2\theta^{q'}.
    \end{aligned}
    $$
    Thus, combining \eqref{Eq: Optimality of -DpH} with A\ref{A: A6} gives
    $$
    \begin{aligned}
        \Lambda_{q'}(\mu)^{q'} & \leq C_0^2\theta^{q'} + \lambda C_0\theta^{q'}\int_{\T^d}\int_0^\infty\left[\int_{\T^d \times \R^d} f(x,m) d(\widetilde{\mu}_t - \mu) + C_0\right]e^{-\lambda t} dtdx \\
        &\leq 4C_0^2\theta^{q'}.
    \end{aligned}
    $$

    For the gradient estimate, we adapt the proof of \cite[Theorem A.1]{lasry1989nonlinear}. We will assume $u$ is smooth as before. Letting $w = |D_xu|^2$, we have
    $$2\lambda w + \nu\left(2|D_{xx}^2u|^2 - \Delta w\right) + D_pH^\theta \cdot D_xw + 2(D_xH^\theta - \theta D_xf) \cdot D_xu = 0.$$
    Letting $x_1 \in \T^d$ be a maximum point for $w$, the maximum principle gives that
    $$
    \begin{aligned}
        \lambda w + \nu|D_{xx}^2u|^2 &\leq  C_0\theta(2 + w^\frac{q}{2} + \theta^{-q'}\Lambda_{q'}(\mu)^{q'})w^{\frac{1}{2}} \\
        &\leq \frac{\lambda}{2}w + \frac{C_0^2\theta^2}{2\lambda}(2 + w^\frac{q}{2} + \theta^{-q'}\Lambda_{q'}(\mu)^{q'})^2
    \end{aligned}
    $$
    at $x_1$. From the Cauchy-Schwarz inequality, we get
    $$\nu|D_{xx}^2u|^2 \geq \frac{\nu}{d}(\Delta u)^2 = \frac{1}{d\nu}(H^\theta(D_xu,\mu) + \lambda u - \theta f)^2.
    $$
    Applying the maximum principle again, we find
    $$\lambda\|u\|_{L^\infty(\T^d)} \leq \|H^\theta(x,0,\mu)\|_{L^\infty(\T^d)} + \theta\|f(x,m)\|_{L^\infty(\T^d)} \leq C\theta.$$
    Combining these results with D\ref{A: Bounds for H} and our bounds on $\Lambda_{q'}(\mu)$ and $\lambda u$, we get that $w^q(x_1) \leq C$. In particular,
    $$\max_x w(x) = w(x_1) \leq C.$$

    For the final estimate, we note that $v(x) \coloneqq u(x) - u(0)$ satisfies
    $$-\nu\Delta v + \int_0^1 D_pH^\theta(x,sD_xv,\mu)ds \cdot D_xv + \lambda v = \theta f(x,m) - H^\theta(x,0,\mu) - \lambda u(0).$$
    By the previous estimates, \cite[Theorem 8.32]{gilbarg1977elliptic} gives a uniform estimate for $[D_xu]_\alpha$ for some $\alpha \in (0,1)$. Thus, since
    $$\nu\Delta v = H^\theta(x,D_xv,\mu) + \lambda(v + u(0)) - \theta f(x,m),$$
    the result follows by \cite[Lemma 6.1]{gilbarg1977elliptic}.
\end{proof}

\begin{remark}
    Observe that for $\theta > 0$ fixed, a nearly identical argument to the one above yields
    $$\Lambda_{q'}(\mu) + \|D_xu\|_\infty \leq C, \qquad\|u - u(0)\|_{C^{2+\alpha(\T^d)}} \leq C(\nu)$$
    uniformly in $\lambda \in (0,1)$.
\end{remark}

With this, we are now ready to prove the main results of this section.

\begin{theorem} \label{Thm: Existence - Stationary}
    Assume A\ref{A: A1}-\ref{A: A6} and D\ref{A: Bounds for H} hold, $d \geq 2$, $\lambda \in (0,1)$, and $\nu \in (0,\frac{1}{d})$. Then there exists at least one solution $(u,m,\mu)$ to \eqref{Eq: Stationary Discounted MFGC}. Moreover, if D\ref{A: U} holds, then there is exactly one solution to \eqref{Eq: Stationary Discounted MFGC}. 
\end{theorem}

\begin{proof}
    The proof is very similar to the one in \cite{graber2025mean,kobeissi2022mean}, and uses the Leray-Schauder fixed point theorem. Define $X$ and $T: X \times [0,1] \to X$ as in Definition \ref{Def: Fixed-Point Map}. Notice that $T(\cdot,\cdot,0) \equiv (0,1)$.
    
    Lemma \ref{Lem: Estimates for Stationary Case} and \ref{Lem: Stationary FP} give estimates for fixed-points. To show that $T$ is continuous, suppose $(\widetilde{u}_k,\widetilde{m}_k,\theta_k)_k$ is a sequence in $X$ converging to some $(\widetilde{u},\widetilde{m},\theta)$.
    Letting $(u_k,m_k) = T(\widetilde{u}_k,\widetilde{m}_k,\theta_k)$, the proof of Lemma \ref{Lem: Estimates for Stationary Case} gives $\|u_k\|_{C^{2+\alpha}(\T^d)} \leq C$ for some $\alpha \in (0,1)$. Hence, there is a subsequence (which we will not rename) that converges to some $u \in C^2(\T^d)$. Adapting the proof of \cite[Lemma 4.4]{kobeissi2022mean}, we find that the measures $\mu_k = (I,-D_pH^\theta(x,D_x\widetilde{u}_k,\mu_k)) \# \eta(\widetilde{m}_k)$ converge to $\mu = (I,-D_pH^\theta(x,D_x\widetilde{\mu},\mu)) \# \eta(\widetilde{m})$. Passing to the limit, we have
    $$-\nu\Delta u + H^\theta(x,D_xu,\mu) + \lambda u = \theta f(x,\eta(\widetilde{m})).$$
    By uniqueness, we get the convergence of the whole sequence.
    
    Similarly, since $(m_k)_k$ is bounded in $L^\infty(\T^d) \cap H^1(\T^d)$, taking a subsequence if necessary, $m_k \to m$ a.e., strongly in $L^r(\T^d)$ for all $r \geq 1$, and weakly in $H^1(\T^d)$. Clearly, we have $m \geq 0$ and $\int_{\T^d} m \, dx = 1$. Now note that for $\varphi \in H^1(\T^d)$, we have
    $$\int_{\T^d} \left[D_x\varphi \cdot (\nu D_xm_k + m_kD_pH(x,D_xu_k,\mu_k)) + \lambda m_k\varphi\right] \, dx = \lambda\int_{\T^d} \varphi dx.$$
    Hence, passing to the limit, we get
    $$\int_{\T^d} \left[D_x\varphi \cdot (\nu D_xm + mD_pH(x,D_xu,\mu)) + \lambda m\varphi\right] \, dx = \lambda\int_{\T^d} \varphi dx.$$
    As before, uniqueness gives convergence of the whole sequence.
    
    That $T$ is compact follows by a similar argument. Thus, by the Leray-Schauder fixed-point theorem, $T$ has a fixed-point, thus concluding the proof of existence.
    The proof of uniqueness is standard (see \cite{graber2025mean,graber2026ergodic,kobeissi2022mean}).
\end{proof}

Now we shift our focus back to the convergence of solutions. To this end, we will require the following estimates for uniform integrability.

\begin{lemma} \label{Lem: Estimates on m for Stationary System}
    Assume A\ref{A: A1}-\ref{A: A6}, C\ref{A: f is semiconcave}, and D\ref{A: Bounds for H}-\ref{A: DpxH = 0 or H is uniformly convex} hold, $d \geq 2$, $\lambda \in (0,1)$, and $\nu \in (0,\frac{1}{d})$. If $(u,m,\mu)$ is a solution to \eqref{Eq: Stationary Discounted MFGC}, then we have
    $$\lambda\int_{\T^d} m\ln{m} \, dx + \nu\int_{\T^d} \frac{|D_xm|^2}{m} dx \leq C$$
    for some constant $C$ independent of $\lambda,\nu$.
\end{lemma}

\begin{proof}
    For $\delta \in (0,1)$, define $\phi_\delta$ to be a smooth function such that:
    \begin{itemize}
        \item $\phi_\delta(\xi) = \ln{\xi}$ for $\xi > \delta$,
        \item $\phi_\delta'(\xi) \geq 0$,
        \item $\xi\phi_\delta''(\xi) + \phi_\delta'(\xi)$ is bounded (uniformly in $\delta$) for $\xi < 1$,
        \item $\phi_\delta(\xi) \downarrow \ln{\xi}$ monotonically as $\delta \to 0$,
        \item $\phi_\delta'(\xi) \uparrow \xi^{-1}$ monotonically as $\delta \to 0$,
        \item and $\phi_\delta''(\xi) \to -\xi^{-2}$ as $\delta \to 0$.
    \end{itemize}
    Multiplying the Fokker-Planck equation by $\phi_\delta(m)$ and integrating, we get
    $$
    \begin{aligned}
        &\lambda \int_{\T^d} (m-1)\phi_\delta(m) dx + \nu\int_{\T^d} |D_xm|^2\phi_\delta'(m) dx = - \int_{\T^d} m\phi_\delta'(m) D_xm \cdot D_pH(x,D_xu,\mu) dx \\
        &= \int_{\T^d} m\bigg((m\phi_\delta''(m) + \phi_\delta'(m))D_xm \cdot D_pH(x,D_xu,\mu) \\
        &\qquad+ m\phi_\delta'(m)\operatorname{Tr}\big(D_{px}^2H(x,D_xu,\mu) + D_{pp}^2H(x,D_xu,\mu)D_{xx}^2u\big)\bigg) dx.
    \end{aligned}
    $$
    By the Monotone Convergence Theorem, we get that
    $$\int_{\T^d} (m-1)\phi_\delta(m) dx \to \int_{\T^d} (m-1)\ln{m} dx$$
    and
    $$\int_{\T^d} |D_xm|^2\phi_\delta'(m) dx \to \int_{\T^d} \frac{|D_xm|^2}{m} dx.$$
    Moreover, by the Dominated Convergence Theorem, we get
    $$
    \begin{aligned}
        &\int_{\T^d} m^2\phi_\delta'(m)\operatorname{Tr}(D_{px}^2H(x,D_xu,\mu) + D_{pp}^2H(x,D_xu,\mu)D_{xx}^2u) dx \\
        &\to \int_{\T^d} m\operatorname{Tr}(D_{px}^2H(x,D_xu,\mu) + D_{pp}^2H(x,D_xu,\mu)D_{xx}^2u) dx
    \end{aligned}
    $$
    and
    $$
        \begin{aligned}
            &\limsup_{\delta\to 0^+} \left|\int_{\T^d} m(m\phi_\delta''(m) + \phi'(m))D_xm \cdot D_pH(x,D_xu,\mu) dx\right| \\
            &\qquad \leq \|m\|_\infty\|D_pH(x,D_xu,\mu)\|_\infty\|D_xm\|_{L^2(\T^d)}\limsup_{\delta\to 0^+} \|m\phi_\delta''(m) + \phi'(m)\|_{L^2(\T^d)} = 0.
        \end{aligned}
    $$
    Since $m > 0$ a.e., Jensen's inequality yields that $\int_{\T^d} \ln{m} \, dx \leq 0$. Thus, the result follows from the fact that $D_{xx}^2u \leq C$, $0 \leq D_{pp}^2H \leq C$, and $|\operatorname{Tr}(D_{px}^2H(x,p,\mu))| \leq C$.
\end{proof}

\begin{theorem} \label{Thm: Vanishing Discount}
    Assume A\ref{A: A1}-\ref{A: A6} and D\ref{A: Bounds for H}-\ref{A: DpxH = 0 or H is uniformly convex} hold, that $d \geq 2$, and that $\nu \in (0,\frac{1}{d})$. Then there is at least one solution to the ergodic problem \eqref{Eq: Ergodic MFGC}.
    Moreover, if D\ref{A: U} also holds, then there is exactly one solution to \eqref{Eq: Ergodic MFGC}, up to adding a constant to $u$.
    
    In particular, under these assumptions,  as $\lambda \to 0$,
    \begin{enumerate}
        \item $u_\lambda - u_\lambda(0) \to u_0$ in $C^2(\T^d)$
        \item $\lambda u_\lambda(0) \to \rho_0$
        \item $m_\lambda \to m_0$ strongly in $L^r(\T^d)$ for all $r \geq 1$ and weakly in $H^1(\T^d)$
        \item $\mu_\lambda \to \mu_0$ in $\mathcal{P}(\T^d \times \R^d)$
    \end{enumerate}
    where $(u_\lambda,m_\lambda,\mu_\lambda)$ is the solution to \eqref{Eq: Stationary Discounted MFGC} for $\lambda \in (0,1)$ and $(u_0,\rho_0,m_0,\mu_0)$ is the solution to \eqref{Eq: Ergodic MFGC} with $u_0(0) = 0$.
\end{theorem}

\begin{proof}
    First, define $v_\lambda = u_\lambda - u_\lambda(0)$. By Lemmas \ref{Lem: Stationary HJ} and \ref{Lem: Estimates for Stationary Case}, we have
    $$\lambda\|u_\lambda\|_{L^\infty(\T^d)} + \|v_\lambda\|_{C^{2+\alpha}(\T^d)} + \Lambda_{q'}(\mu_\lambda) \leq C$$
    uniformly in $\lambda$ for some $\alpha \in (0,1)$. By the Arzela-Ascoli Theorem and by the tightness of measures, there exist $v \in C^2(\T^d)$, $\rho \in \R$, and $\mu \in \mathcal{P}(\T^d \times \R^d)$ such that, passing to a subsequence if necessary, we get $v_\lambda \to v$ in $C^2(\T^d)$, $\lambda u_\lambda(0) \to \rho$ in $\R$, and $\mu_\lambda \to \mu$ in $\mathcal{P}(\T^d \times \R^d)$ as $\lambda \to 0$.
    Furthermore, by Lemma \ref{Lem: Stationary FP}, we have
    $$\|m_\lambda\|_{H^1(\T^d)} + \|m_\lambda\|_{L^\infty(\T^d)} \leq C.$$
    Thus, there exists $m \in H^1(\T^d) \cap L^\infty(\T^d)$ so that, up to a subsequence, $m_\lambda \to m$ strongly in $L^r(\T^d)$ for all $r \geq 1$ and weakly in $H^1(\T^d)$. Passing to the limit, we get that $v$ is a classical solution to $-\nu\Delta v + H(x,D_xv,\mu) + \rho = f(x,m)$, $m$ is a solution to $\nu\Delta m + \operatorname{div} (mD_pH(x,D_xv,\mu)) = 0$ in the weak sense, and $\mu = (I,-D_pH(\cdot,D_xu,\mu)) \# m$.
    
    What remains is to prove the uniqueness of solutions to \eqref{Eq: Ergodic MFGC}. The uniqueness of solutions to the Hamilton-Jacobi equation (up to adding a constant) for fixed $m,\mu$ follows from \cite[Theorem 5]{gomes2002stochastic}, and the uniqueness of solutions (in the sense of distributions) to the Fokker-Planck equation for fixed $u,\mu$ follows from Lemma \ref{Lem: Estimates for Equations in Divergence Form} (see also \cite[Corollary 2]{gomes2002stochastic} or \cite[Section 4.1]{bogachev2022fokker}).
    The proof of uniqueness for solutions to \eqref{Eq: Ergodic MFGC} therefore follows by standard arguments (see \cite{graber2025mean,graber2026ergodic,kobeissi2022mean}).
\end{proof}

\begin{remark} \label{Rmk: Holder Regularity of m}
    By \cite[Section 1.7]{bogachev2022fokker}, we obtain bounds for the solution $m_0$ to the ergodic Fokker-Planck equation in $C^\beta(\T^d)$.
\end{remark}

\begin{theorem} \label{Thm: Vanishing Viscosity}
    Assume A\ref{A: A1}-\ref{A: A6}, C\ref{A: f is semiconcave}, and D\ref{A: Bounds for H}-\ref{A: DpxH = 0 or H is uniformly convex} hold, that $d \geq 2$, and that $\lambda \in (0,1)$. Then there exists at least one solution $(u_0,m_0,\mu_0)$ to \eqref{Eq: Stationary First-Order MFGC}. Moreover, if
    \begin{equation} \label{Eq: Unique m}
        \int_{\T^d} (f(x,m_1) - f(x,m_2)) d(m_1 - m_2) > 0, \qquad \forall m_1,m_2 \in \mathcal{P}(\T^d) \text{ with } m_1 \neq m_2
    \end{equation}
    then there is exactly one solution $(u_0,m_0,\mu_0)$ to \eqref{Eq: Stationary First-Order MFGC}, and as $\nu \to 0$, we have that $u_\nu \to u_0$ in $C^0(\T^d) \cap W^1_r(\T^d)$ for all $r > 1$, $m_\nu \to m_0$ weakly in $L^1(\T^d)$, and $\mu_\nu \to \mu_0$ in $\mathcal{P}(\T^d \times \R^d)$, where $(u_\nu,m_\nu,\mu_\nu)$ is the solution to \eqref{Eq: Stationary Discounted MFGC} for each $\nu \in (0,\frac{1}{d})$. In the first-order case, the HJ equation is understood in the viscosity sense and the FP equation is understood in the sense of distributions.
\end{theorem}

\begin{proof}
    By Lemmas \ref{Lem: Stationary HJ} and \ref{Lem: Estimates for Stationary Case}, we have
    $$\|u_\nu\|_{C^1(\T^d)} + \Lambda_{q'}(\mu_\nu) \leq C \qquad \text{and} \qquad D_{xx}^2u_\nu \leq C$$
    uniformly in $\nu$. Hence, passing to a subsequence if necessary, we that as $\nu \to 0$, $u_\nu$ converges uniformly to some $u \in C^0(\T^d) \cap W^1_\infty(\T^d)$ and $\mu_\nu$ converges to some $\mu \in \mathcal{P}(\T^d \times \R^d)$. Moreover, by \cite[Theorem 3.3.3]{cannarsa2004semiconcave} and our uniform semiconcavity estimate, we conclude that $D_xu_\nu \to D_xu$ a.e. Combining this with our gradient bound implies $u_\nu \to u$ in $W^1_r(\T^d)$ for all $r \geq 1$.
    Similarly, by Lemma \ref{Lem: Estimates on m for Stationary System}, we get that
    $$\lambda\int_{\T^d} m_\nu\ln(m_\nu) \, dx \leq C,$$
    and hence $(m_\nu)_\nu$ is uniformly integrable. By the Dunford-Pettis Theorem, there exists some $m \in L^1(\T^d)$ such that, up to a subsequence, $m_\nu \to m$ weakly in $L^1$.

    Recalling that for all $\varphi \in H^1(\T^d)$, we have
    $$\nu\int_{\T^d} D_x\varphi \cdot D_xu_\nu dx + \int_{\T^d} \varphi (H(x,D_xu_\nu,\mu_\nu) + \lambda u_\nu) dx = \int_{\T^d} \varphi f(x,m_\nu) dx,$$
    passing to the limit (and taking a subsequence if necessary), we get that
    $$\int_{\T^d} \varphi (H(x,D_xu,\mu) + \lambda u) dx = \int_{\T^d} \varphi f(x,m) dx$$
    for all $\varphi \in H^1(\T^d)$.
    Thus, $u$ is a weak solution to $H(x,D_xu,\mu) + \lambda u = f(x,m)$. Moreover, $u$ is also a viscosity solution by stability.
    
    Similarly, for each $\varphi \in C^2(\T^d)$, we have
    $$
    \begin{aligned}
        \lambda\int_{\T^d} \varphi dx &= \int_{\T^d} [\lambda\varphi - \nu\Delta\varphi + D_x\varphi \cdot D_pH(x,D_xu_\nu,\mu_\nu)] dm_\nu \\
        &= \int_{\T^d} [\lambda\varphi - \nu\Delta\varphi + D_x\varphi \cdot D_pH(x,D_xu,\mu)] dm_\nu \\
        &\qquad+ \int_{\T^d} D_x\varphi \cdot (D_pH(x,D_xu_\nu,\mu_\nu) - D_pH(x,D_xu,\mu)) dm_\nu.
    \end{aligned}
    $$
    As $\nu \to 0$, we have that
    $$\int_{\T^d} [\lambda\varphi - \nu\Delta\varphi + D_x\varphi \cdot D_pH(x,D_xu,\mu)] dm_\nu \to \int_{\T^d} [\lambda\varphi + D_x\varphi \cdot D_pH(x,D_xu,\mu)] dm,$$
    and passing to a subsequence if necessary, we have $D_pH(x,D_xu_\nu,\mu_\nu) \to D_pH(x,D_xu,\mu)$ a.e. By the uniform bounded on $D_pH(x,D_xu_\nu,\mu_\nu)$ and the uniform integrability of $m_\nu$, for every $\varepsilon > 0$, there exists $M_\varepsilon > 0$ such that
    $$\int_{m_\nu \geq M_\varepsilon} |D_x\varphi \cdot (D_pH(x,D_xu_\nu,\mu_\nu) - D_pH(x,D_xu,\mu))| dm_\nu \leq \varepsilon.$$
    Defining $m_{\nu,\varepsilon} = m_\nu\chi_{m_\nu \leq M_\varepsilon}$, we conclude by the Dominated Convergence Theorem that for each $\varepsilon > 0$,
    $$\lim_{\nu \to 0} \int_{\T^d} D_x\varphi \cdot (D_pH(x,D_xu_\nu,\mu_\nu) - D_pH(x,D_xu,\mu)) m_{\nu,\varepsilon} dx = 0.$$
    Therefore,
    $$
    \begin{aligned}
        &\limsup_{\nu \to 0} \left|\int_{\T^d} D_x\varphi \cdot (D_pH(x,D_xu_\nu,\mu_\nu) - D_pH(x,D_xu,\mu)) dm_\nu\right| \\
        &\leq \limsup_{\nu \to 0} \left|\int_{\T^d} D_x\varphi \cdot (D_pH(x,D_xu_\nu,\mu_\nu) - D_pH(x,D_xu,\mu)) d(m_\nu - m_{\nu,\varepsilon})\right| \\
        &\qquad + \limsup_{\nu \to 0} \left|\int_{\T^d} D_x\varphi \cdot (D_pH(x,D_xu_\nu,\mu_\nu) - D_pH(x,D_xu,\mu)) dm_{\nu,\varepsilon}\right| \\
        &\leq \varepsilon.
    \end{aligned}
    $$
    Letting $\varepsilon \to 0$ we conclude that
    $$\int_{\T^d} D_x\varphi \cdot (D_pH(x,D_xu_\nu,\mu_\nu) - D_pH(x,D_xu,\mu)) dm_\nu \to 0.$$
    Thus, $\lambda m - \operatorname{div} (mD_pH(x,D_xu,\mu)) = \lambda$ in the sense of distributions.
    By the continuity of $D_pH$, we have $\mu = (I,-D_pH(\cdot,D_xu,\mu)) \# m$.
    
    As before, we conclude the proof by showing uniqueness. 
    By standard arguments (see \cite{graber2025mean,graber2026ergodic,kobeissi2022mean}), if $(u^1,m^1,\mu^1)$ and $(u^2,m^2,\mu^2)$ are two solutions, then we find
    $$
    \begin{aligned}
        0 \leq& \int_\Omega (f(x,m^1) - f(x,m^2)) d(m^1 - m^2) \\
        =&\int_\Omega m^1(L(x,\alpha^{\mu^2},\mu^2) - L(x,\alpha^{\mu^1},\mu^1) + D_\alpha L(x,\alpha^{\mu^2},\mu^2) \cdot (\alpha^{\mu^1}-\alpha^{\mu^2})) dx \\
        &+ \int_\Omega m^2(L(x,\alpha^{\mu^1},\mu^1) - L(x,\alpha^{\mu^2},\mu^2) + D_\alpha L(x,\alpha^{\mu^1},\mu^1) \cdot (\alpha^{\mu^2}-\alpha^{\mu^1})) dx.
    \end{aligned}
    $$
    Since $L$ is strictly convex,
    \begin{equation}
        L(x,\alpha_1,\mu^i) - L(x,\alpha_2,\mu^i) + D_\alpha L(x,\alpha_1,\mu^i) \cdot (\alpha_2 - \alpha_1) \leq 0
        \label{Eq: Convexity of L}
    \end{equation}
    with equality holding if and only if $\alpha_1 = \alpha_2$. Hence,
    $$
    \begin{aligned}
        0 \leq& \int_\Omega \bigg(m^1(L(x,\alpha^{\mu^1},\mu^2) - L(x,\alpha^{\mu^1},\mu^1)) + m^2(L(x,\alpha^{\mu^2},\mu^1) - L(x,\alpha^{\mu^2},\mu^2))\bigg)dx \\
        =& -\int_{\Omega \times \R^d} (L(x,\alpha,\mu^1) - L(x,\alpha,\mu^2)) d(\mu^1 - \mu^2)(x,\alpha).
    \end{aligned}
    $$
    By A\ref{A: A3}, this gives
    $$\int_\Omega \bigg(m^1(L(x,\alpha^{\mu^1},\mu^2) - L(x,\alpha^{\mu^1},\mu^1)) + m^2(L(x,\alpha^{\mu^2},\mu^1) - L(x,\alpha^{\mu^2},\mu^2))\bigg)dx = 0.$$
    By the condition for equality for \eqref{Eq: Convexity of L}, we get $\alpha^{\mu^1} = \alpha^{\mu^2}$. Moreover, by \eqref{Eq: Unique m}, we have $m^1 = m^2$ and hence $\mu^1 = (I,\alpha^{\mu^i}) \# m^i = \mu^2$. By the uniqueness of solutions to the Hamilton-Jacobi equation (e.g., see \cite[Theorem 1.14]{tran2021hamilton}), we find $u^1 = u^2$.
\end{proof}

We observe that the uniqueness of $\alpha_0 = -D_pH(x,D_xu_0,\mu_0)$ and $f(x,m_0)$ can be obtained under weaker monotonicity assumptions. To conclude this section, we analyze the rate of the convergence found in Theorem \ref{Thm: Vanishing Viscosity}.

\begin{theorem} \label{Thm: Convergence Rate for Vanishing Viscosity}
    Assume A\ref{A: A1}-\ref{A: A6}, C\ref{A: f is semiconcave}, and D\ref{A: Bounds for H}-\ref{A: U} hold, that $d \geq 2$, and that $\lambda \in (0,1)$. Then for all $\nu \in [0,\frac{1}{d})$,
    $$\int_{\T^d} (f(x,m_\nu) - f(x,m_0))(m_\nu - m_0) dx \leq C\sqrt{\nu}$$
    and if D\ref{A: L is Uniformly Convex} also holds, then
    $$\int_{\T^d} (m_\nu + m_0)|\alpha_\nu - \alpha_0|^2 dx \leq C\sqrt{\nu}$$
    where $(u_\nu,m_\nu,\mu_\nu)$ is a solution to \eqref{Eq: Stationary Discounted MFGC} and $\alpha_\nu = -D_pH(x,D_xu_\nu,\mu_\nu)$.
\end{theorem}

\begin{proof}
    First, note that
    $$\lambda(u_\nu - u_0) - \nu\Delta u_\nu + H(x,D_xu_\nu,\mu_\nu) - H(x,D_xu_0,\mu_0) = f(x,m_\nu) - f(x,m_0).$$
    Multiplying by $m_\nu - m_0$, integrating, and recalling that
    $$\lambda(m_\nu - m_0) -\nu\Delta m_\nu - \operatorname{div}(m_\nu D_pH(x,D_xu_\nu,\mu_\nu) - m_0D_pH(x,D_xu_0,\mu_0)) = 0$$
    in the sense of distributions, we obtain
    $$
    \begin{aligned}
        &\int_{\T^d} (f(x,m_\nu) - f(x,m_0))(m_\nu - m_0) dx \\
        &= \int_{\T^d} m_\nu\left(H(x,D_xu_\nu,\mu_\nu) - H(x,D_xu_0,\mu_0) + (D_xu_0 - D_xu_\nu) \cdot D_pH(x,D_xu_\nu,\mu_\nu)\right) dx \\
        &\qquad + \int_{\T^d} m_0\left(H(x,D_xu_0,\mu_0) - H(x,D_xu_\nu,\mu_\nu) + (D_xu_\nu - D_xu_0) \cdot D_pH(x,D_xu_0,\mu_0)\right) dx \\
        &\qquad + \nu\int_{\T^d} (m_0\Delta u_\nu + D_xm_\nu \cdot D_x u_0) dx.
    \end{aligned}
    $$
    Letting $\alpha_\nu \coloneqq -D_pH(x,D_xu_\nu,\mu_\nu)$, we recall that
    $$L(x,\alpha_\nu,\mu_\nu) = D_xu_\nu \cdot D_pH(x,D_xu_\nu,\mu_\nu) - H(x,D_xu_\nu,\mu_\nu), \qquad D_xu_\nu = -D_\alpha L(x,\alpha_\nu,\mu_\nu).$$
    Thus, using A\ref{A: A3}, we obtain
    $$
    \begin{aligned}
        &\int_{\T^d} (f(x,m_\nu) - f(x,m_0))(m_\nu - m_0) dx \\
        &+ \int_{\T^d} m_\nu\left(L(x,\alpha_\nu,\mu_0) - L(x,\alpha_0,\mu_0) + (\alpha_0 - \alpha_\nu) \cdot D_\alpha L(x,\alpha_0,\mu_0)\right) dx \\
        &+ \int_{\T^d} m_0\left(L(x,\alpha_0,\mu_\nu) - L(x,\alpha_\nu,\mu_\nu) + (\alpha_\nu - \alpha_0) \cdot D_\alpha L(x,\alpha_\nu,\mu_\nu)\right) dx \\
        &\qquad \leq \nu\int_{\T^d} (m_0\Delta u_\nu + D_xm_\nu \cdot D_x u_0) dx.
    \end{aligned}
    $$
    What remains is to estimate the righthand side. First, recall that $\Delta u_\nu \leq C$ (see Lemma \ref{Lem: Stationary HJ}). Thus,
    $$\nu\int_{\T^d} m_0\Delta u_\nu dx \leq C\nu.$$
    Now recall that for $\nu \in (0,1)$, we have $m_\nu \in H^1(\T^d)$. Therefore, applying Lemma \ref{Lem: Estimates on m for Stationary System}, we obtain
    $$\nu\int_{\T^d} D_xm_\nu \cdot D_xu_0 dx\leq \nu\left(\int_{\T^d} \frac{|D_xm_\nu|^2}{m_\nu}dx\right)^\frac{1}{2} \left(\int_{\T^d} |D_xu_0|^2m_\nu dx\right)^\frac{1}{2} \leq C\sqrt{\nu}.$$
    Thus, the first inequality follows from the fact that $L$ is strictly convex in $\alpha$, and if D\ref{A: L is Uniformly Convex} holds, then we have the second inequality.
\end{proof}

%-----%-----%-----%
%-----%-----%-----%
%-----%-----%-----%
%-----%-----%-----%
\section{Long-Time Average Behavior}
\label{Sec: Long-Time Average}
%-----%-----%-----%
%-----%-----%-----%
%-----%-----%-----%
%-----%-----%-----%

In this section, we investigate the long-time average behavior of solutions to the zero-discount MFGC system \eqref{Eq: MFGC with No Discount Factor}. We prove that for a certain class of example Hamiltonians, solutions to \eqref{Eq: MFGC with No Discount Factor} converge to the solution to the ergodic problem \eqref{Eq: Ergodic MFGC}, and we prove some weaker convergence results for more general Hamiltonians.

\begin{lemma} \label{Lem: Estimates for uT}
    Assume A\ref{A: A1}-\ref{A: A6} hold. Let $T \geq 1$ and let $(u^T,m^T,\mu^T)$ be a solution to \eqref{Eq: MFGC with No Discount Factor} in the sense of distributions. Then there exists some constant $C > 0$ independent of $T$ such that
    $$\int_0^T \Lambda_{q'}(\mu^T)^{q'} dt + \|u^T\|_\infty \leq CT.$$
    Moreover, if D\ref{A: No x-dependence} holds, then
    $$\|D_xu^T\|_\infty + \Lambda_{q'}(\mu^T) \leq C.$$
    Alternatively, if D\ref{A: Example Hamiltonian}, then we have
    $$\Lambda_{q'}(\mu^T)^{q'} + \|D_xu^T\|_\infty \leq C, \qquad D_{xx}^2u^T \leq C.$$
\end{lemma}

\begin{proof}
    From \cite{kobeissi2022mean}, we get
    $$\int_0^T \Lambda_{q'}(\mu^T)^{q'} dt + \|u^T\|_\infty \leq CT.$$
    Now suppose D\ref{A: No x-dependence} holds. As in the proof of Lemma \ref{Lem: General C1 Estimate}, define
    $$\varphi(v) = \exp(\exp(-1-R^{-1}v)) \qquad \text{for } |v| \leq R$$
    $$w(t,x) = \varphi\left(T^{-1}u^T(t,x) + \frac{(1 + C_0)}{T}\int_t^T \Lambda_{q'}(\mu^T)^{q'} ds\right)|D_xu^T(t,x)|^2$$
    for a sufficiently large constant $R > 0$. Then we get
    $$
    \begin{aligned}
        &-\partial_tw - \nu\Delta w + D_xw \cdot D_pH + \frac{2\nu\varphi'}{T\varphi}D_xw \cdot D_xu^T + \frac{2\nu\varphi}{T}|D_{xx}^2u^T|^2 \\
        &= \frac{\varphi'}{T\varphi}w\left[f + D_xu^T \cdot D_pH - H + (1+C_0)\Lambda_{q'}(\mu^T)^{q'}\right] - \nu\frac{\varphi''\varphi - 2(\varphi')^2}{T^2\varphi^3}w^2 \\
        &\leq T^{-1}\left(\varepsilon - \frac{1}{C}\right)\left(w^{1+\frac{q}{2}} + w\Lambda_{q'}(\mu^T)^{q'}\right) + C_{R,\varepsilon}T^{-1}
    \end{aligned}
    $$
    and hence the maximum principle yields that
    $$\|w\|_\infty \leq \|w(T,\cdot)\|_\infty + \int_0^T CT^{-1} dt \leq C.$$
    Therefore, $|D_xu^T| \leq Cw^\frac{1}{2} \leq C$. By \cite[Lemma 4.1]{kobeissi2022mean}, this gives $\Lambda_{q'}(\mu^T) \leq C$.

    Finally, suppose $H(p,\mu) = \psi(\mu)|p + \varphi(\mu)|^2 + V(\mu)$ and
    $$\underset{m}{\sup}(\|f(\cdot,m)\|_{W^2_\infty(\T^d)} + \|g(\cdot,m)\|_{C^{3+\beta}(\T^d)}) \leq C_0.$$
    We conclude by proving our semiconcavity estimate, which will give us a uniform estimate for the gradient. We will again assume that our solutions are smooth, as otherwise we would perform a convolution (we note as in \cite{kobeissi2022mean}, that we can apply \cite[Theorem 6.48]{lieberman1996second} to get that $D_{xx}^2u^T$ and $D_{xxx}^3u^T$ are continuous).

    Fix $\xi \in \R^d$ with $|\xi| = 1$. Note that
    $$-\partial_t\partial_{\xi\xi}^2 u^T - \nu\Delta\partial_{\xi\xi}^2 u^T + 2\psi\left(|D_x\partial_\xi u^T\right|^2 + (D_xu^T + \varphi) \cdot D_x\partial_{\xi\xi}^2 u^T) = \partial_{\xi\xi}^2 f.$$
    By the maximum principle, at a point $(t,x)$ such that $\partial_{\xi\xi}^2 u^T = \max \partial_{\xi\xi}^2 u^T_+$, either $t = T$ and so $\partial_{\xi\xi}^2 u^T \leq \|g(\cdot,m^T(T))\|_{C^2(\T^d)} \leq C_0$, or else
    $$|\partial_{\xi\xi}^2 u^T_+|^2 \leq |D_x\partial_\xi u^T|^2 \leq \frac{1}{2\psi}\|\partial_{\xi\xi}^2f\|_\infty \leq C.$$
    Taking the supremum over all such $\xi$, we conclude the proof.
\end{proof}

In this section, to simplify notation, we will define the following:
\begin{equation}
    b^T \coloneqq -D_pH(x,D_xu^T,\mu^T), \qquad b = -D_pH(x,D_xu,\mu)
\end{equation}
and
\begin{equation}
    U^T \coloneqq u^T - (u + \rho(T-t)), \qquad M^T \coloneqq m^T - m.
\end{equation}

\begin{lemma} \label{Lem: Energy Estimates for Long-Time Problem}
    Assume A\ref{A: A1}-\ref{A: A6}, D\ref{A: Bounds for H}-\ref{A: DpxH = 0 or H is uniformly convex}, and D\ref{A: L is Uniformly Convex} hold, $d \geq 2$, and $\nu \in (0,\frac{1}{d})$. Let $T \geq 1$, let $(u^T,m^T,\mu^T)$ solve \eqref{Eq: MFGC with No Discount Factor}, and let $(u,\rho,m,\mu)$ solve \eqref{Eq: Ergodic MFGC} with $\int_{\T^d} u \, dx = 0$. Then for $0 \leq t_1 \leq t_2 \leq T$
    $$
    \begin{aligned}
        &\int_{t_1}^{t_2}\int_{\T^d} \left[\frac{m^T + m}{C}\left|b^T - b\right|^2 + (f(x,m^T)- f(x,m))M^T\right] dxdt \\
        &\hspace{3cm} + \int_{t_1}^{t_2}\int_{\T^d} (L(\alpha,\mu^T) - L(\alpha,\mu))d(\mu^T - \mu) dt \leq -\left[\int_{\T^d} M^TU^T dx\right]_{t_1}^{t_2}.
    \end{aligned}
    $$
\end{lemma}

\begin{proof}
    For ease of presentation, we will not write $x$ as an argument for $f$. Now note that
    $$
    \begin{aligned}
        &\left[\int_{\T^d} M^TU^T dx\right]_{t_1}^{t_2} \\
        &= \int_{t_1}^{t_2}\int_{\T^d} \big[-D_xU^T \cdot (\nu D_xm^T - m^Tb^T) \\
        &\qquad+ (m^T - m)(-\nu\Delta u^T + H(D_xu^T,\mu^T) - f(m^T) + \rho)\big] dxdt \\
        &= \int_{t_1}^{t_2}\int_{\T^d} \big[m^T(H(D_xu^T,\mu^T) - H(D_xu,\mu) + D_xU^T \cdot b^T) \\
        &\qquad+ m(H(D_xu,\mu) - H(D_xu^T,\mu^T) - D_xU^T \cdot b) + M^T(f(m) - f(m^T))\big] dxdt \\
        &= \int_{t_1}^{t_2}\int_{\T^d} \big[m^T(L(b,\mu) - L(b^T,\mu^T) + (b^T - b) \cdot D_\alpha L(b,\mu)) \\
        &\qquad+ m(L(b^T,\mu^T) - L(b,\mu) + (b - b^T) \cdot D_\alpha L(b^T,\mu^T)) + (m^T - m)(f(m) - f(m^T))\big] dxdt \\
    \end{aligned}
    $$
    By D\ref{A: L is Uniformly Convex}, we get
    $$L(b^T,\mu^T) - L(b,\mu) - (b^T - b) \cdot D_\alpha L(b,\mu) \geq \frac{1}{C}|b^T - b|^2 + L(b^T,\mu^T) - L(b^T,\mu)$$
    and
    $$L(b,\mu) - L(b^T,\mu^T) - (b - b^T) \cdot D_\alpha L(b^T,\mu^T) \geq \frac{1}{C}|b^T - b|^2 + L(b,\mu) - L(b,\mu^T).$$
    Thus, we have
    $$
    \begin{aligned}
        &\int_{t_1}^{t_2}\int_{\T^d} \left[\frac{m^T + m}{C}\left|b^T - b\right|^2 + (f(m^T) - f(m))M^T\right] dxdt  \\
        &\hspace{3cm} + \int_{t_1}^{t_2}\int_{\T^d} (L(\alpha,\mu^T) - L(\alpha,\mu))d(\mu^T - \mu) dt\leq -\left[\int_{\T^d} M^TU^T dx\right]_{t_1}^{t_2}.
    \end{aligned}
    $$
\end{proof}

\begin{lemma} \label{Lem: Long-Time Convergence of b}
    Assume that A\ref{A: A1}-\ref{A: A6}, D\ref{A: Bounds for H}-\ref{A: L is Uniformly Convex}, and either D\ref{A: No x-dependence} or D\ref{A: Example Hamiltonian} hold, that $d \geq 2$, and that $\nu \in (0,\frac{1}{d})$. Now let $T \geq 1$, let $(u^T,m^T,\mu^T)$ be a solution to \eqref{Eq: MFGC with No Discount Factor}, and let $(u,\rho,m,\mu)$ be a solution to \eqref{Eq: Ergodic MFGC} with $\int_{\T^d} u \, dx = 0$. Then
    \begin{equation} \label{Eq: Energy Estimate for Long-Time Problem}
        \begin{aligned}
            &\int_0^T\int_{\T^d} \left[\frac{m^T + m}{C}\left|b^T - b\right|^2 + (f(x,m^T)- f(x,m))M^T\right] dxdt \\
            &\hspace{3.5cm}+ \int_0^T\int_{\T^d} (L(\alpha,\mu^T) - L(\alpha,\mu))d(\mu^T - \mu) dt \leq C
        \end{aligned}
    \end{equation}
    and hence
    \begin{equation} \label{Eq: b^T converges in L^2(m)}
        \lim_{T \to \infty} \int_0^1\int_{\T^d} (m^T(Tt,\cdot) + m)\left|b^T(Tt,\cdot) - b\right|^2 dxdt = 0.
    \end{equation}
\end{lemma}

\begin{proof}
    By Lemma \ref{Lem: Energy Estimates for Long-Time Problem}, the result follows from the fact that
    $$\int_{\T^d} (m^T(T) - m)(u - u^T(T)) dx \leq 2(\|g\|_\infty + \|u\|_\infty) \leq C$$
    and
    $$\int_{\T^d} (m^T(0) - m)(u^T(0) - u) dx \leq 2\|u\|_\infty + \operatorname{osc}u^T(0) \leq 2\left(\|u\|_\infty + \|D_xu^T\|_\infty\right) \leq C.$$
\end{proof}

To obtain stronger convergence results, we observe that
\begin{equation}
    \begin{cases}
        -\partial_tU^T - \nu\Delta U^T + H(D_xu^T,\mu^T) - H(D_xu,\mu) = f(m^T) - f(m) \\
        \partial_tM^T - \nu\Delta M^T + \operatorname{div} (M^Tb^T) = \operatorname{div} (m(b - b^T)) \\
        M^T(0,x) = m_0(x) - m(x), \qquad U^T(T,x) = g(x,m^T(T)) - u(x)
    \end{cases}
\end{equation}
in the sense of distributions. Now we recall the following consequence of \cite[Lemma 7.6]{cardaliaguet2013long}:

\begin{lemma} \label{Lem: Estimates for Equations in Divergence Form}
    Let $V \in L^\infty([0,T] \times \T^d;\R^d)$, $F_1 \in L^2([0,T] \times \T^d)$, and $F_2 \in C^0([0,T];L^2(\T^d))$ with $\int_{\T^d} F_2(t,x)dx = 0$ for all $t$. If $\eta \in L^1([0,T] \times \T^d)$ is a solution (in the sense of distributions) to
    $$\partial_t\eta - \nu\Delta\eta - \operatorname{div} (\eta V) = \operatorname{div} (F_1) + F_2, \qquad \text{in } (0,T) \times \T^d$$
    with $\int_{\T^d} \eta(t,x) dx = 0$ for all $t \in (0,T)$, then for every $0 \leq t_1 \leq t_2 \leq T$,
    $$\|\eta(t_2)\|_{L^2(\T^d)} \leq C\left[\|\eta(t_1)\|_{L^2(\T^d)}e^{-\omega(t_2-t_1)} + \int_{t_1}^{t_2}e^{-\omega(t_2-s)}\|F_2(s)\|_{L^2(\T^d)} ds + \|F_1\|_{L^2([t_1,t_2] \times \T^d)}\right]$$
    for some constants $C,\omega > 0$ depending only on $\nu$ and $\|V\|_\infty$.
\end{lemma}

\begin{lemma} \label{Lem: First Estimates on MT}
    Assume A\ref{A: A1}-\ref{A: A6}, D\ref{A: Bounds for H}-\ref{A: L is Uniformly Convex}, and either D\ref{A: No x-dependence} or D\ref{A: Example Hamiltonian} hold, that $d \geq 2$, and that $\nu \in (0,\frac{1}{d})$. Let $T \geq 1$, let $(u^T,m^T,\mu^T)$ be a solution to \eqref{Eq: MFGC with No Discount Factor}, and let $(u,\rho,m,\mu)$ be the solution to \eqref{Eq: Ergodic MFGC} with $\int_{\T^d} u \, dx = 0$. Then there exist constants $C,\omega > 0$ (depending only on $\nu$ and the constants in the assumptions) such that for every $0 \leq t_0 \leq t \leq T$,
    \begin{equation} \label{Eq: L2 Estimate for MT}
        \|M^T(t)\|_{L^2(\T^d)} \leq C\left(\|M^T(t_0)\|_{L^2(\T^d)}e^{-\omega(t-t_0)} + \|m(b^T - b)\|_{L^2([t_0,t] \times \T^d)}\right)
    \end{equation}
    and hence
    $$\sup_t\|M^T(t)\|_{L^2(\T^d)} \leq C.$$
\end{lemma}

\begin{proof}
    The first inequality is a straightforward consequence of Lemma \ref{Lem: Estimates for Equations in Divergence Form}. For the second, we take $t_0 = 0$ and use Lemma \ref{Lem: Long-Time Convergence of b} and Remark \ref{Rmk: Holder Regularity of m} to obtain
    $$\|m(b^T - b)\|_{L^2([0,T] \times \T^d)}^2 \leq \|m\|_{L^\infty(\T^d)}\int_0^T\int_{\T^d} m|b^T - b|^2 dx dt \leq C.$$
\end{proof}

\begin{theorem}
    Assume A\ref{A: A1}-\ref{A: A6}, D\ref{A: Bounds for H}-\ref{A: L is Uniformly Convex}, and either D\ref{A: No x-dependence} or D\ref{A: Example Hamiltonian} hold, that $d \geq 2$, and that $\nu \in (0,\frac{1}{d})$. Now let $(u^T,m^T,\mu^T)$ be solutions to \eqref{Eq: MFGC with No Discount Factor} and $(u,\rho,m,\mu)$ be the solution to \eqref{Eq: Ergodic MFGC} with $\int_{\T^d} u \, dx = 0$. As $T \to \infty$, we have the following:
    \begin{enumerate}
        \item the map $b^T(Tt,\cdot)$ converges to $b$ in $L^r([0,1] \times \T^d)$ and $L^r(m)$ for all $r \geq 1$.
        \item the map $m^T(Tt,\cdot)$ converges to $m$ weakly in $L^2([0,1] \times \T^d)$.
        \item the map $\mu^T(Tt)$ converges to $\mu$ in $\mathcal{P}([0,1] \times \T^d \times \R^d)$ in the sense that for any continuous function $\phi$, we have $\int_{[0,1] \times \T^d \times \R^d} \phi d(\mu^T - \mu) \to 0$.
    \end{enumerate}
\end{theorem}

\begin{proof}
    By Lemma \ref{Lem: Long-Time Convergence of b}, $b^T(Tt,\cdot)$ converges to $b$ in $L^2(m)$, and the uniform bound on $b^T$ gives convergence in $L^r(m)$ for all $r \geq 1$. Moreover, taking a subsequence if necessary, we get $b^T(Tt,\cdot) \to b$ $m$-a.e. By the Harnack Inequality for weak solutions (see \cite[Theorem 8.20]{gilbarg1977elliptic}), we have $m > 0$ a.e., and hence this implies that $b^T(Tt,\cdot) \to b$ a.e. in $[0,1] \times \T^d$. By the uniform bound on $|b^T - b|$, we conclude that $b^T \to b$ in $L^r([0,1] \times \T^d)$ (and the uniqueness of the $L^r$ limit gives convergence of the entire sequence).
    
    By Lemma \ref{Lem: First Estimates on MT}, we have that $m^T(Tt,\cdot)$ is bounded in $L^2([0,1] \times \T^d)$. Since this is a separable, reflexive Banach space, we can extract a subsequence converging weakly to some $\widetilde{m}$. Now recall that for all $\phi \in C^\infty_c([0,1) \times \T^d)$, we have
    $$\int_0^1\int_{\T^d} \left(-\frac{1}{T}\partial_t\phi - \nu\Delta\phi - D_x\phi \cdot b^T\right)m^T(Tt,\cdot) dxdt = \frac{1}{T}\int_{\T^d}m_0\phi(0) dx.$$
    Passing to the limit, we have
    $$\int_0^1\int_{\T^d} \left(-\nu\Delta\phi + D_x\phi \cdot b\right)\widetilde{m} dxdt = 0$$
    for all $\phi \in C^\infty_c([0,1) \times \T^d)$. Fixing $\phi \in C^2(\T^d)$ and $t_0 \in (0,1)$, we construct a sequence of smooth functions $\varphi_\delta(t,x) = \phi_\delta(x)\eta_\delta(t)$ as follows. Choose $\phi_\delta \in C^\infty(\T^d)$ such that $\phi_\delta \to \phi$ in $H^2(\T^d)$ and choose $\eta_\delta(t) = \delta^{-1}\eta\left(\frac{t-t_0}{\delta}\right)$ for $\delta < \min\{t_0,1-t_0\}$, where
    $$\eta(s) = 
    \begin{cases}
        e^{-\frac{1}{1-s^2}}, &|s| < 1 \\
        0, &|s| \geq 1
    \end{cases}
    $$
    Then for a.e. $t \in (0,1)$,
    $$\int_{\T^d} (-\nu\Delta\phi - D_x\phi \cdot b)\widetilde{m}dx = \lim_{\delta \to 0} \int_0^1 \int_{\T^d} (-\nu\Delta\phi_\delta - D_x\phi_\delta \cdot b)\widetilde{m}dxdt = 0.$$
    As this holds for all $\phi \in C^2(\T^d)$, we conclude that $\widetilde{m}(t,\cdot) = m$ for a.e. $t$. By uniqueness of the limit, we conclude that $m^T \to m$ weakly in $L^2([0,1] \times \T^d)$.
    
    For the convergence of $\mu^T(Tt)$, fix $\phi \in C^0([0,1] \times \T^d \times \R^d)$. Then we observe that
    $$
    \begin{aligned}
        &\limsup_{T \to \infty} \left|\int_{[0,1] \times \T^d \times \R^d} \phi(t,x,\alpha) d(\mu^T - \mu)(t,x,\alpha)\right| \\
        &\leq \limsup_{T \to \infty} \int_0^1\int_{\T^d} |\phi(t,x,-D_pH(D_xu^T,\mu^T)) - \phi(t,x,-D_pH(D_xu,\mu))| dm^T dt \\
        &\qquad + \limsup_{T \to \infty} \left|\int_0^1\int_{\T^d} \phi(t,x,-D_pH(D_xu,\mu)) d(m^T - m) dt\right| \\
        &\leq \sup_T\|m^T\|_{L^2}\limsup_{T \to \infty}\|\phi(t,x,-D_pH(D_xu^T,\mu^T)) - \phi(t,x,-D_pH(D_xu,\mu))\|_{L^2} \\
        &= 0,
    \end{aligned}
    $$
    where the last step follows from Lemma \ref{Lem: First Estimates on MT} and the dominated convergence theorem.
\end{proof}

If we can obtain some notion of equicontinuity of $m^T(Tt,\cdot)$ in time, then this result can be strengthened naturally to achieve much stronger convergence. However, this has proved to be rather difficult in general due to the rescaling of the time horizon. Instead, we will focus our attention on a specific case and adapt an argument from \cite{cardaliaguet2013long} to obtain stronger convergence. In particular, we will consider the following case:
\begin{assumptionD} \label{special case}
    $$
\begin{cases}
    H(p,\mu) = c|p + \varphi(\mu)|^2 + V(\mu) \\
    L(\alpha,\mu) = \frac{1}{4c}|\alpha|^2 + \alpha \cdot \varphi(\mu) - V(\mu)
\end{cases}
$$
where $V,\varphi$ are continuous functions such that
\begin{enumerate}
    \item $|V(\mu)| + |\varphi(\mu)| \leq C_0(1 + \Lambda_2(\mu)^2)$;
    \item $(\varphi(\mu_1) - \varphi(\mu_2)) \cdot \int_{\T^d \times \R^d} \widetilde{\alpha} d(\mu_1 - \mu_2) \geq 0$, which we note is equivalent to A\ref{A: A3};
    \item $\varphi(\mu_k) \to \varphi(\mu)$ whenever $(\varphi(\mu_k) - \varphi(\mu)) \cdot \int_{\T^d \times \R^d} \widetilde{\alpha} d(\mu_k - \mu) \to 0$.
\end{enumerate}
\end{assumptionD}
We observe that the linear demand model found in \cite{chan2015bertrand} is such an example.

First, we will prove an analogue of \cite[Lemma 3.6]{cardaliaguet2013long}.

\begin{lemma} \label{Lem: Long-Time Convergence of m}
    Assume A\ref{A: A1}-\ref{A: A6}, D\ref{A: Bounds for H}-\ref{A: L is Uniformly Convex}, and D\ref{A: Example Hamiltonian}-\ref{special case} hold, that $d \geq 2$, and that $\nu \in (0,\frac{1}{d})$. Let $T \geq 1$, let $(u^T,m^T,\mu^T)$ be a solution to \eqref{Eq: MFGC with No Discount Factor}, and let $(u,\rho,m,\mu)$ be the solution to \eqref{Eq: Ergodic MFGC} with $\int_{\T^d} u \, dx = 0$. Then for every $\delta \in (0,\frac{1}{2})$, we have
    \begin{equation} \label{Eq: Long-Time Convergence of m}
        \int_{\delta}^{1-\delta}\int_{\T^d} |b^T(Tt) - b|^2dm(x) dt \to 0 \qquad \text{and} \qquad \sup_{t \in [\delta,1-\delta]}\|M^T(Tt)\|_{L^2(\T^d)} \to 0.
    \end{equation}
    In particular, for each $t \in (0,1)$, $m^T(Tt,\cdot) \to m$ in $L^2(\T^d)$.
\end{lemma}

\begin{proof}
    By \eqref{Eq: Energy Estimate for Long-Time Problem}, we have that
    $$\int_0^1 \left((\varphi(\mu^T(Tt)) - \varphi(\mu)) \cdot \int_{\T^d} (b^T(Tt)m^T(Tt) - bm) dx\right) dt \to 0,$$
    and D\ref{special case} says that the integrand is non-negative, so it converges in $L^1(0,1)$.
    Taking a subsequence if necessary, by D\ref{special case} we also get $\varphi(\mu^T(Tt)) \to \varphi(\mu)$ for a.e. $t \in [0,1]$.
    Also by Lemma \ref{Lem: Long-Time Convergence of b}, we get
    $$\int_0^T\int_{\T^d} (m + m^T)||b^T - b|^2 dxdt \leq C.$$
    Thus, we can choose some $t_1 \in [0,(\delta/2)]$ and $t_2 \in [(1-\delta/2),1]$ such that $\varphi(\mu^T(Tt_1)) \to \varphi(\mu)$, $\varphi(\mu^T(Tt_2)) \to \varphi(\mu)$, and
    $$\int_{\T^d} (m + m^T(t_1))|b^T(t_1) - b|^2 dx \to 0, \qquad \int_{\T^d} (m + m^T(t_2))|b^T(t_2) - b|^2 dx \to 0.$$
    
    Applying \eqref{Eq: Energy Estimate for Long-Time Problem}, we obtain
    $$\int_{t_1}^{t_2}\int_{\T^d} (m + m^T)|b^T - b|^2 dxdt \leq \int_{\T^d} M^T(t_1)U^T(t_1) dx - \int_{\T^d} M^T(t_2)U^T(t_2) dx.$$
    Letting $V^T(t,x) = U^T(t,x) - \int_{\T^d} U^T(t,y) dy$, since $m \geq \frac{1}{C} > 0$, we obtain
    $$
    \begin{aligned}
        \int_{\T^d} M^T(t_2)U^T(t_2) dx &\leq \int_{\T^d} (m^T(t_2) + m)|V^T(t_2)| dx \\
        &\leq \|V^T(t_2)\|_{L^2(\T^d)}\|m^T(t_2) + m\|_{L^2(\T^d)} \\
        &\leq C\|D_xU^T(t_2)\|_{L^2(\T^d)}\|m^T(t_2) + m\|_{L^2(\T^d)} \\
        &\leq C\left(\|b^T(t_2) - b\|_{L^2(\T^d)} + |\varphi(\mu^T(Tt_2) - \varphi(\mu)|\right)\|m^T(t_2) + m\|_{L^2(\T^d)} \\
        &\leq C\left(\|b^T(t_2) - b\|_{L^2(m)} + |\varphi(\mu^T(Tt_2) - \varphi(\mu)|\right) \\
        &\to 0
    \end{aligned}
    $$
    as $T \to \infty$ (for a fixed $\delta > 0$). Applying a similar argument to $\int_{\T^d} M^T(t_1)U^T(t_1) dx$, we obtain
    $$\int_{\delta T}^{(1-\delta)T}\int_{\T^d} |b^T - b|^2dm(x) dt \leq \int_{t_1}^{t_2}\int_{\T^d} |b^T - b|^2dm(x) dt \to 0.$$
    Finally, combining this with \eqref{Eq: L2 Estimate for MT} gives
    $$
    \begin{aligned}
        &\limsup_{T \to \infty} \sup_{t \in [\delta, 1-\delta]}\|M^T(Tt)\|_{L^2(\T^d)} \\
        &\leq C\limsup_{T \to \infty} \left(\|M^T(t_1)\|_{L^2(\T^d)}e^{-\omega(\delta/2)T} + \|m\|_{L^\infty(\T^d)}^\frac{1}{2}\left(\int_{t_1}^{t_2}\int_{\T^d} |b^T - b|^2dm(x) ds\right)^\frac{1}{2}\right) \\
        &= 0,
    \end{aligned}
    $$
    thus completing the proof.
\end{proof}

With this, we are now ready to prove our final convergence result.

\begin{theorem} \label{Thm: Long-Time Average for Example}
    Assume A\ref{A: A1}-\ref{A: A6}, D\ref{A: Bounds for H}-\ref{A: Lipschitz Assumption on DpH Inverse}, and D\ref{A: Example Hamiltonian}-\ref{special case} hold, that $d \geq 2$, and that $\nu \in (0,\frac{1}{d})$. Let $T \geq 1$, let $(u^T,m^T,\mu^T)$ be a solution to \eqref{Eq: MFGC with No Discount Factor}, and let $(u,\rho,m,\mu)$ be the solution to \eqref{Eq: Ergodic MFGC} with $\int_{\T^d} u \, dx = 0$. Then
    \begin{enumerate}
        \item the map $m^T(Tt,\cdot)$ converges to $m$ strongly in $L^2([0,1] \times \T^d)$;
        \item the map $\mu^T(Tt)$ converges to $\mu$ in $L^r(0,1;\mathcal{P}(\T^d \times \R^d))$ for $r \geq 1$;
        \item the map $u^T(Tt,\cdot) - \int_{\T^d} u^T(Tt,y) dy$ converges to $u$ in $W^{0,1}_r([0,1] \times \T^d)$ for $r \geq 1$;
        \item the map $u^T(Tt,\cdot)/T$ converges to $(1-t)\rho$ in $L^r([0,1] \times \T^d$) for $r \geq 1$.
    \end{enumerate}
\end{theorem}

\begin{proof}
    The convergence of $m^T$ is straightforward: applying Lemmas \ref{Lem: First Estimates on MT} and \ref{Lem: Long-Time Convergence of m}, for each $\delta \in (0,\frac{1}{2})$, we have
    $$
    \begin{aligned}
        \|M^T(Tt)\|_{L^2([0,1] \times \T^d)}^2 &\leq \|M^T(Tt)\|_{L^2([\delta T,(1 - \delta)T]  \times \T^d)}^2 + \|M^T(Tt)\|_{L^2(([0,1] \setminus[\delta T,(1 - \delta)T])  \times \T^d)}^2\\
        &\leq \|M^T(Tt)\|_{L^2([\delta T,(1 - \delta)T]  \times \T^d)}^2 + 2C\delta \\
        &\to 2C\delta.
    \end{aligned}
    $$
    as $T \to \infty$, and letting $\delta \to 0$ gives the desired convergence.

    Turning now to $\mu^T$, we note that for $t \in [0,1]$, we have
    $$
    \begin{aligned}
        W_1(\mu^T(Tt),\mu) &= \sup_{Lip(\varphi) \leq 1} \bigg[\int_{\T^d} (\varphi(x,b^T(Tt)) - \varphi(x,b))dm \\
        &\qquad\qquad + \int_{\T^d} (\varphi(x,b^T(Tt)) - \varphi(0,0)) d(m^T(Tt) - m)\bigg] \\
        &\leq \int_{\T^d} |b^T(Tt)-b|dm + \int_{\T^d} (|x| + |b^T(Tt)|)|m^T(Tt) - m| dx
    \end{aligned}
    $$
    and so
    $$\int_0^1 W_1(\mu^T(Tt),\mu) dt \leq \int_0^1 \int_{\T^d}\left(|b^T(Tt)-b|m + C|m^T(Tt) - m|\right) dxdt \to 0.$$
    This shows that $\mu^T(Tt) \to \mu$ in $L^1(0,1;\mathcal{P}(\T^d \times \R^d))$, and since $W_1(\mu^T(Tt),\mu)$ is uniformly bounded, it follows that $\mu^T(Tt) \to \mu$ in $L^r(0,1;\mathcal{P}(\T^d \times \R^d))$ for all $r \geq 1$.
    
    By assumption D\ref{A: Lipschitz Assumption on DpH Inverse}, we get that (up to a subsequence)
    $$
    \begin{aligned}
        |D_xu^T(Tt,\cdot) - D_xu| \leq C\left(|b^T(Tt,\cdot) - b| + |b + D_pH(D_xu,\mu^T(Tt))|\right) \to 0
    \end{aligned}
    $$
    pointwise a.e. By our uniform bounds on $|D_xu^T(Tt,\cdot) - D_xu|$, we get that $D_xu^T(Tt,\cdot) \to D_xu$ in $L^r([0,1] \times \T^d)$ for all $r \geq 1$. Thus, by Poincar\'{e}'s inequality, we get that $u^T(Tt,\cdot) - \int_{\T^d} u^T(Tt,y) dy$ converges to $u$ in $L^r([0,1] \times \T^d)$ for all $r \geq 1$.

    Integrating the Hamilton-Jacobi equation satisfied by $u^T$, we find that for $t \in [0,1]$,
    $$
    \begin{aligned}
        \frac{1}{T}\int_{\T^d} u^T(Tt,\cdot)dx &= \frac{1}{T}\int_{\T^d} g(\cdot,m^T(T)) dx \\
        &\qquad + \int_t^1\int_{\T^d} \left(f(m^T(Ts)) -H(D_xu^T(Ts),\mu^T(Ts))\right) dxds \\
        &\to \int_t^1\int_{\T^d} \left(f(m) -H(D_xu,\mu)\right) dxds \\
        &= \rho(1-t)
    \end{aligned}
    $$
    uniformly in $t$ as $T \to \infty$. Thus, by the boundedness of $u^T(Tt,\cdot) - \int_{\T^d} u^T(Tt,y) dy$, we have that
    $$
    \begin{aligned}
        \left\|\frac{1}{T}u^T(Tt,x) - \rho(1-t)\right\|_{L^r([0,1] \times \T^d)} &\leq \frac{1}{T}\left\|u^T(Tt,x) - \int_{\T^d} u^T(Tt,y)dy\right\|_{L^r([0,1] \times \T^d)} \\
        &\qquad + \left\|\frac{1}{T}\int_{\T^d} u^T(Tt,y)dy - \rho(1-t)\right\|_{L^r([0,1] \times \T^d)} \\
        &\to 0
    \end{aligned}
    $$
    as $T \to \infty$ for all $r \geq 1$.
\end{proof}

%-----%-----%-----%
%-----%-----%-----%
%-----%-----%-----%
%-----%-----%-----%
\section*{Funding \& Conflicts of Interest}
%-----%-----%-----%
%-----%-----%-----%
%-----%-----%-----%
%-----%-----%-----%

The authors are grateful to be supported by National Science Foundation through NSF Grant DMS-2045027. The authors also declare that they have no conflicts of interest.

%-----%-----%-----%
%-----%-----%-----%
%-----%-----%-----%
%-----%-----%-----%
\section*{Acknowledgements}
%-----%-----%-----%
%-----%-----%-----%
%-----%-----%-----%
%-----%-----%-----%

We would like to thank Son Tu for his technical assistance with the well-posedness of the stationary Hamilton-Jacobi and Fokker-Planck equations.

%-----%-----%-----%
%-----%-----%-----%
%-----%-----%-----%
%-----%-----%-----%
\section*{Disclosure of AI Use}
%-----%-----%-----%
%-----%-----%-----%
%-----%-----%-----%
%-----%-----%-----%

The authors disclose that Claude Opus 5.5 was used to check for errors after the article was completed.

\bibliographystyle{amsplain}
\bibliography{refs}

@article{kobeissi2022classical,
  title={On classical solutions to the mean field game system of controls},
  author={Kobeissi, Ziad},
  journal={Communications in Partial Differential Equations},
  volume={47},
  number={3},
  pages={453--488},
  year={2022},
  publisher={Taylor \& Francis}
}

@article{kobeissi2022mean,
  title={Mean field games with monotonous interactions through the law of states and controls of the agents},
  author={Kobeissi, Ziad},
  journal={Nonlinear Differential Equations and Applications NoDEA},
  volume={29},
  number={5},
  pages={52},
  year={2022},
  publisher={Springer}
}

@article{bongini2024mean,
  title={Mean Field Games of controls with Dirichlet boundary conditions},
  author={Bongini, Mattia and Salvarani, Francesco},
  journal={ESAIM: Control, Optimisation and Calculus of Variations},
  volume={30},
  pages={32},
  year={2024},
  publisher={EDP Sciences}
}

@book{ladyzhenskaia1968linear,
  title={Linear and quasi-linear equations of parabolic type},
  author={Ladyzhenskaia, Olga Aleksandrovna and Solonnikov, Vsevolod Alekseevich and Ural'tseva, Nina N},
  volume={23},
  year={1968},
  publisher={American Mathematical Soc.}
}

@article{lasry2007mean,
  title={Mean field games},
  author={Lasry, Jean-Michel and Lions, Pierre-Louis},
  journal={Japanese journal of mathematics},
  volume={2},
  number={1},
  pages={229--260},
  year={2007},
  publisher={Springer}
}

@article{Huang2006Large,
author = {Minyi Huang and Roland P. Malham{\'e} and Peter E. Caines},
title = {Large population stochastic dynamic games: closed-loop McKean-Vlasov systems and the Nash certainty equivalence principle},
volume = {6},
journal = {Communications in Information \& Systems},
number = {3},
publisher = {International Press of Boston},
pages = {221 -- 252},
year = {2006},
}

@article{gomes2016extended,
  title={Extended deterministic mean-field games},
  author={Gomes, Diogo A and Voskanyan, Vardan K},
  journal={SIAM Journal on Control and Optimization},
  volume={54},
  number={2},
  pages={1030--1055},
  year={2016},
  publisher={SIAM}
}

@article{cardaliaguet2018mean,
  title={Mean field game of controls and an application to trade crowding},
  author={Cardaliaguet, Pierre and Lehalle, Charles-Albert},
  journal={Mathematics and Financial Economics},
  volume={12},
  pages={335--363},
  year={2018},
  publisher={Springer}
}

@article{cardaliaguet2013long,
  title={Long time average of mean field games with a nonlocal coupling},
  author={Cardaliaguet, Pierre and Lasry, Jean-Michel and Lions, Pierre-Louis and Porretta, Alessio},
  journal={SIAM Journal on Control and Optimization},
  volume={51},
  number={5},
  pages={3558--3591},
  year={2013},
  publisher={SIAM}
}

@article{bo2025mean,
  title={Mean field game of controls with state reflections: Existence and limit theory},
  author={Bo, Lijun and Wang, Jingfei and Yu, Xiang},
  journal={arXiv preprint arXiv:2503.03253},
  year={2025}
}

@article{graber2025mean,
  title={Mean Field Games of Controls with Boundary Conditions \& Invariance Constraints},
  author={Graber, P Jameson and Rosengartner, Kyle},
  journal={arXiv preprint arXiv:2508.21642},
  year={2025}
}

@article{graber2026ergodic,
  title={Ergodic Mean Field Games of Controls with State Constraints},
  author={Graber, P Jameson and Rosengartner, Kyle},
  journal={arXiv preprint arXiv:2604.07550},
  year={2026}
}

@article{gomes2014existence,
  title={On the existence of classical solutions for stationary extended mean field games},
  author={Gomes, Diogo A and Patrizi, Stefania and Voskanyan, Vardan},
  journal={Nonlinear Analysis: Theory, Methods \& Applications},
  volume={99},
  pages={49--79},
  year={2014},
  publisher={Elsevier}
}

@article{lasry1989nonlinear,
      title={Nonlinear elliptic equations with singular boundary conditions and stochastic control with state constraints: 1. The Model problem},
      author={Lasry, Jean-Michel and Lions, Pierre-Louis},
      journal={Mathematische Annalen},
      volume={283},
      number={4},
      pages={583--630},
      year={1989},
      publisher={Springer}
}

@book{gilbarg1977elliptic,
  title={Elliptic partial differential equations of second order},
  author={Gilbarg, David and Trudinger, Neil S},
  volume={224},
  number={2},
  year={1977},
  publisher={Springer}
}

@article{graber2021note,
  title={A note on mean field games of controls with state constraints: existence of mild solutions},
  author={Graber, Jameson and Mayorga, Sergio},
  journal={arXiv preprint arXiv:2109.11655},
  year={2021}
}

@article{bardi2021convergence,
  title={Convergence of some mean field games systems to aggregation and flocking models},
  author={Bardi, Martino and Cardaliaguet, Pierre},
  journal={Nonlinear Analysis},
  volume={204},
  pages={112199},
  year={2021},
  publisher={Elsevier}
}

@article{bertucci2019some,
  title={Some remarks on mean field games},
  author={Bertucci, Charles and Lasry, Jean-Michel and Lions, Pierre-Louis},
  journal={Communications in Partial Differential Equations},
  volume={44},
  number={3},
  pages={205--227},
  year={2019},
  publisher={Taylor \& Francis}
}

@book{lieberman1996second,
  title={Second order parabolic differential equations},
  author={Lieberman, Gary M},
  year={1996},
  publisher={World scientific}
}

@book{cannarsa2004semiconcave,
  title={Semiconcave Functions, Hamilton—Jacobi Equations, and Optimal Control},
  author={Cannarsa, Piermarco and Sinestrari, Carlo},
  year={2004},
  publisher={Springer}
}

@article{cardaliaguet2012long,
  title={Long time average of mean field games.},
  author={Cardaliaguet, Pierre and Lasry, Jean-Michel and Lions, Pierre-Louis and Porretta, Alessio},
  journal={Networks \& Heterogeneous Media},
  volume={7},
  number={2},
  year={2012}
}

@book{tran2021hamilton,
  title={Hamilton--Jacobi equations: theory and applications},
  author={Tran, Hung V},
  volume={213},
  year={2021},
  publisher={American Mathematical Soc.}
}

@article{gomes2002stochastic,
  title={A stochastic analogue of Aubry-Mather theory},
  author={Gomes, Diogo Aguiar},
  journal={Nonlinearity},
  volume={15},
  number={3},
  pages={581--603},
  year={2002}
}

@article{chan2015bertrand,
  title={Bertrand and Cournot mean field games},
  author={Chan, Patrick and Sircar, Ronnie},
  journal={Applied Mathematics \& Optimization},
  volume={71},
  number={3},
  pages={533--569},
  year={2015},
  publisher={Springer}
}

@article{cardaliaguet2015mean,
  title={Mean field games systems of first order},
  author={Cardaliaguet, Pierre and Graber, P Jameson},
  journal={ESAIM: Control, Optimisation and Calculus of Variations},
  volume={21},
  number={3},
  pages={690--722},
  year={2015},
  publisher={EDP Sciences}
}

@article{degond2017mean,
  title={Mean-field games and model predictive control},
  author={Degond, Pierre and Herty, Michael and Liu, Jian-Guo},
  journal={Communications in Mathematical Sciences},
  volume={15},
  number={5},
  pages={1403--1422},
  year={2017},
  publisher={International Press of Boston}
}

@article{graber2026mean,
  title={Mean field games of controls with fractional laplacian},
  author={Graber, P Jameson and Matter, Elizabeth and Bolanos, Jesus Ruiz},
  journal={Nonlinear Differential Equations and Applications NoDEA},
  volume={33},
  number={3},
  pages={64},
  year={2026},
  publisher={Springer}
}

@article{graber2024new,
  title={New uniqueness results for a mean field game of controls},
  author={Graber, Jameson and Matter, Elizabeth},
  journal={arXiv preprint arXiv:2410.13950},
  year={2024}
}

@article{chan2017fracking,
  title={Fracking, Renewables, and Mean Field Games},
  author={Chan, Patrick and Sircar, Ronnie},
  journal={SIAM Review},
  url={https://doi.org/10.1137/15M1031424},
  volume={59},
  number={3},
  pages={588--615},
  year={2017},
  publisher={SIAM}
}

@article{graber2018existence,
  title={{Existence and uniqueness of solutions for Bertrand and Cournot mean field games}},
  author={Graber, P Jameson and Bensoussan, Alain},
  journal={Applied Mathematics \& Optimization},
  url={https://doi.org/10.1007/s00245-016-9366-0},
  volume={77},
  number={1},
  pages={47--71},
  year={2018},
  publisher={Springer}
}

@inproceedings{graber2018variational,
  title={Variational mean field games for market competition},
  author={Graber, P Jameson and Mouzouni, Charafeddine},
  booktitle={PDE Models for Multi-Agent Phenomena},
  url={https://arxiv.org/abs/1707.07853},
  series={Springer INdAM Series},
  year={2018},
  volume={28},
}

@article{graber2020mean,
  title={On mean field games models for exhaustible commodities trade},
  author={Graber, Philip Jameson and Mouzouni, Charafeddine},
  journal={ESAIM: Control, Optimisation and Calculus of Variations},
  volume={26},
  pages={11},
  year={2020},
  publisher={EDP Sciences}
}

@article{graber2023master,
  title={Master equation for Cournot mean field games of control with absorption},
  author={Graber, P Jameson and Sircar, Ronnie},
  journal={Journal of Differential Equations},
  volume={343},
  pages={816--909},
  year={2023},
  publisher={Elsevier}
}

@article{camilli2025learning,
  title={Learning equilibria in Cournot mean field games of controls},
  author={Camilli, Fabio and Lauri{\`e}re, Mathieu and Tang, Qing},
  journal={SIAM Journal on Control and Optimization},
  volume={63},
  number={2},
  pages={1407--1431},
  year={2025},
  publisher={SIAM}
}

@article{bonnans2019schauder,
  title={Schauder Estimates for a Class of Potential Mean Field Games of Controls},
  author={Bonnans, J Fr{\'e}d{\'e}ric and Hadikhanloo, Saeed and Pfeiffer, Laurent},
  journal={Applied Mathematics \& Optimization},
  url={https://doi.org/10.1007/s00245-019-09592-z},
  pages={1--34},
  year={2019},
  publisher={Springer}
}

@article{tang2025convergence,
  title={The convergence rate of vanishing viscosity approximations for mean field games},
  author={Tang, Wenpin and Zhang, Yuming Paul},
  journal={SIAM Journal on Mathematical Analysis},
  volume={57},
  number={3},
  pages={3217--3254},
  year={2025},
  publisher={SIAM}
}

@article{gomes2019selection,
  title={The selection problem for some first-order stationary mean-field games},
  author={Gomes, Diogo A and Mitake, Hiroyoshi and Terai, Kengo},
  journal={arXiv preprint arXiv:1908.06485},
  year={2019}
}

@article{terai2026remarks,
  title={Remarks on selection problems for first-order discounted mean field games},
  author={Terai, Kengo},
  journal={arXiv preprint arXiv:2606.03580},
  year={2026}
}

@article{cardaliaguet2013KAM,
  title={Long time average of first order mean field games and weak KAM theory},
  author={Cardaliaguet, Pierre},
  journal={Dynamic Games and Applications},
  volume={3},
  number={4},
  pages={473--488},
  year={2013},
  publisher={Springer}
}

@article{bardi2024long,
  title={Long-time behavior of deterministic mean field games with nonmonotone interactions},
  author={Bardi, Martino and Kouhkouh, Hicham},
  journal={SIAM Journal on Mathematical Analysis},
  volume={56},
  number={4},
  pages={5079--5098},
  year={2024},
  publisher={SIAM}
}

@article{cirant2024long,
  title={Long time behavior and stabilization for displacement monotone mean field games},
  author={Cirant, Marco and M{\'e}sz{\'a}ros, Alp{\'a}r R},
  journal={arXiv preprint arXiv:2412.14903},
  year={2024}
}

@article{cardaliaguet2019long,
  title={Long time behavior of the master equation in mean field game theory},
  author={Cardaliaguet, Pierre and Porretta, Alessio},
  journal={Analysis \& PDE},
  volume={12},
  number={6},
  pages={1397--1453},
  year={2019},
  publisher={Mathematical Sciences Publishers}
}

@article{cirant2025uniqueness,
  title={Uniqueness of solutions to MFG systems with large discount},
  author={Cirant, Marco and Continelli, Elisa},
  journal={arXiv preprint arXiv:2510.09280},
  year={2025}
}

@article{ciampa2026rate,
  title={On the rate of the vanishing viscosity approximation for Mean Field Games with nonlocal coupling},
  author={Ciampa, Gennaro and Goffi, Alessandro},
  journal={arXiv preprint arXiv:2607.07526},
  year={2026}
}

@article{iturriaga2025discretization,
  title={Discretization and Vanishing Discount Problems for First-order Mean Field Games},
  author={Iturriaga, Renato and Mendico, Cristian and Wang, Kaizhi and Xu, Yuchen},
  journal={arXiv preprint arXiv:2509.14541},
  year={2025}
}

@article{piccoli2013transport,
  title={Transport equation with nonlocal velocity in Wasserstein spaces: convergence of numerical schemes},
  author={Piccoli, Benedetto and Rossi, Francesco},
  journal={Acta applicandae mathematicae},
  volume={124},
  number={1},
  pages={73--105},
  year={2013},
  publisher={Springer}
}

@book{bogachev2022fokker,
  title={Fokker--Planck--Kolmogorov Equations},
  author={Bogachev, Vladimir I and Krylov, Nicolai V and R{\"o}ckner, Michael and Shaposhnikov, Stanislav V},
  volume={207},
  year={2022},
  publisher={American Mathematical Society}
}

@article{cannarsa2023weak,
  title={Weak kam approach to first-order mean field games with state constraints},
  author={Cannarsa, Piermarco and Cheng, Wei and Mendico, Cristian and Wang, Kaizhi},
  journal={Journal of Dynamics and Differential Equations},
  volume={35},
  number={2},
  pages={1885--1916},
  year={2023},
  publisher={Springer}
}

@article{jackson2026unconditional,
  title={Unconditional well-posedness of the master equation for monotone mean field games of controls},
  author={Jackson, Joe and M{\'e}sz{\'a}ros, Alp{\'a}r R},
  journal={arXiv preprint arXiv:2601.22845},
  year={2026}
}

@article{mou2022mfgc,
  title={Mean field games of controls: Propagation of monotonicities},
  author={Mou, Chenchen and Zhang, Jianfeng},
  journal={Probability, Uncertainty and Quantitative Risk},
  volume={7},
  number={3},
  year={2022}
}

@article{li2026master,
  title={Master Equations for Mean Field Game of Controls: A Unification of Weak Solution Notions},
  author={Li, Mengzhen and Liu, Xintian and Mou, Chenchen and Wu, Zhen},
  journal={arXiv preprint arXiv:2609.16524},
  year={2026}
}

\end{document}